\documentclass[smallextended,referee,envcountsect,]{svjour3}
\usepackage{amsmath}
\usepackage{lscape}
\usepackage{pdflscape}
\usepackage{subcaption}
\usepackage{amssymb}
\smartqed
\usepackage{graphicx}
\usepackage[ruled,vlined]{algorithm}
\usepackage{algorithmicx}
\usepackage{algpseudocode}
\journalname{JOTA}

\begin{document}

\title{Yuille-Poggio's Flow and Global Minimizer of Polynomials through Convexification by Heat Evolution}


\author{Qiao Wang}

\institute{Qiao Wang \at School of Information Science and Engineering, Southeast University \\
              Nanjing, 210096, China\\
              \and
              School of Economics and Management, Southeast University \\
              Nanjing, 210096, China\\
qiaowang@seu.edu.cn
}

\date{Received: date / Accepted: date}

\maketitle

\begin{abstract}
This study examines the convexification version of the backward differential flow algorithm for the global minimization of polynomials, introduced by O. Arikan \textit{et al} in \cite{ABK}. It investigates why this approach might fail with high-degree polynomials yet succeeds with quartic polynomials. We employ the heat evolution method for convexification combined with Gaussian filtering, which acts as a cumulative form of Steklov's regularization. In this context, we apply the fingerprint theory from computer vision. Originally developed by A.L. Yuille and T. Poggio in the 1980s for computer vision, the fingerprint theory, particularly the fingerprint trajectory equation, is used to illustrate the scaling (temporal) evolution of minimizers. In the case of general polynomials, our research has led to the creation of the Yuille-Poggio flow and a broader interpretation of the fingerprint concepts, in particular we establish the condition both sufficient and necessary  for the convexified backward differential flow algorithms to successfully achieve global minimization. For quartic polynomials, our analysis not only reflects the results of O. Arikan et al. \cite{ABK} but also presents a significantly simpler version of Newton's method that can always globally minimize quartic polynomials without  convexification.

\end{abstract}
\keywords{convex optimization \and Yuille-Poggio's flow \and fingerprint theory\and polynomial optimization \and heat equation \and multiscale Gaussian filter \and computer vision \and Euclidean algorithm}
\subclass{35Q90 \and 46N10 \and 35Q93  \and 90C26}


\section{Background and Motivations}

Global optimization of real polynomials is an important non-convex optimization problem (cf. \cite{Las} and references there in), and produces many excellent theories in past decades.  Among them, N. Z. Shor \cite{Shor-1} first transformed univariate polynomial optimization to convex problem through quadratic optimization in 1987, which can offer approximate solution to this global optimization. After that, N.Z. Shor further studied its relationship with Hilbert's 17th problem \cite{Shor-2}. Also in 1987, V. N. Nefnov \cite{Nef} proposed an algorithm by computing the roots of algebraic equation for finding the minimizer.

In 2014, J. Zhu, S. Zhao and G. Liu \cite{Zhu-2} proposed a backward differential flow formulation, comes from Kuhn-Tucker equation of constrained optimization, to find out the global minimizer of polynomials. They consider the problem for sufficient smooth function $p(x)$,
\begin{equation}\label{zhu-11}
\begin{split}
&\min p(x)\\
&s.t.\ \ x\in \mathcal D:=\{x\in \mathbb R^n|\ \|x\|<a\}
\end{split}
\end{equation}
by introducing a set
\begin{equation}
G=\{\rho>0|\ [\nabla^2 p(x)+\rho I ]>0,\ \forall x\in \mathcal D\},
\end{equation}
and an initial pair $(\widehat \rho,\widehat x)\in G\times\mathcal D$ satisfying
\begin{equation}
\nabla p(\widehat x)+\widehat \rho\widehat x=0.
\end{equation}
Then they proved that the back differential flow $\widehat x(\rho)$,  defined near $\widehat \rho$,
\begin{equation}
\begin{split}
\frac{\,d\widehat x}{\,d\rho}+&[\nabla^2 p(\widehat x)+\rho I]^{-1}\widehat x=0,\\
 &\widehat x(\widehat \rho)=\widehat x
\end{split}
\end{equation}
will lead to the solution of  \eqref{zhu-11}.

 The above work is under the condition that all global minimizers of this polynomial occur only in a known ball, thus the unconstrained optimization problem may be reduced to a constrained optimization problem. However, O.~Arikan, R.S.~Burachik and C.Y.~Kaya \cite {ABK-1} pointed out in 2015 that the method in \cite{Zhu-2} may not converge to global minimizer by a counter-example of quartic polynomial
\begin{equation}
p(x)=x^4-8x^3-18x^2+56x.
\end{equation}
 Moreover, they \cite{ABK} introduced a Steklov regularization and trajectory approach for optimizing univariate polynomials in 2019. Subsequently, in 2020, R.S.~Burachik and C.Y.~Kaya \cite{BK} extended this method to handle multivariable cases. In these studies, quartic polynomial optimization serves as intriguing toy examples. Additionally, polynomials of degree six and higher often do not reach global minimizers, as demonstrated through various examples and counter-examples in \cite{ABK}.


 Actually, the Steklov regularization \cite{ABK}
 \begin{equation}
 \mu(x,t)=\frac{1}{2t}\int_{x-t}^{x+t}f(\tau)\,d\tau
 \end{equation}
 is a low-pass filter, in the viewpoint of signal processing, since we may  write
  \begin{equation}\label{eqn:steklov}
 \mu(x,t)=\frac{1}{2t}\int_{x-t}^{x+t}f(\tau)\,d\tau=f(x)*1_{[-t,t]}(x)
 \end{equation}
 where
 \begin{equation}
 1_{[-t,t]}(x)=\begin{cases}\frac{1}{2t}, &{x\in [-t,t]}\\
 0,&{x\notin [-t,t]}\\
 \end{cases}
 \end{equation}
  from which one may obtain $\mu_x(x,t),\mu_{xt}(x,t)$ and $\mu_{xx}(x,t)$ (where the subscript means partial derivative). However,  $\mu_t(x,t)$ is not explicitly in this regime, although we can represent a differential equation
  \begin{equation}
  2\mu_t+t\mu_{tt}=t\mu_x.
  \end{equation}
Obviously it brings some inconvenience in analyzing the evolution of local minimizers.
 Therefore, we require an approach which can  balance between the simple differential equation and filters, as well as  offer convexification for polynomials.
Fortunately, the heat conduct equation
\begin{equation}\label{a0}
\frac{\partial p}{\partial t}=\frac{1}{2}\frac{\partial^2 p}{\partial x^2}
\end{equation}
with initial condition
\begin{equation}\label{a0-ic}
    p(x,0)=p(x)
\end{equation}
is a nice framework to implement the convexification and the similar optimization algorithm. In addition, the analysis for evolution of all critical points becomes more analytically tractable.

\begin{remark}
It should be pointed that the initial value problem of heat equation \eqref{a0} is equivalent to Gaussian filter, which will be explained in Subsection \ref{subsection-x}. But on the other hand, the accumulation of Steklov regularization will lead to Gaussian distribution, since that
\begin{equation}
\underbrace{1_{[-t,t]}*1_{[-t,t]}*\cdots * 1_{[-t,t]}}_n\to  
\frac{1}{t\sqrt{2\pi n}} \exp\left(-\frac{x^2}{2nt^2}\right)
\end{equation}
for $n$ large enough. Thus replacing Steklov regularity with heat evolution, i.e., the Gaussian filtering, is very natural in this paper.
\end{remark}

Our interest in this paper is to explore the method of optimizing the  even degree monic polynomial
\begin{equation}\label{a1}
    \min_x p(x)=x^n+\sum_{j=1}^n c_jx^{n-j}.
\end{equation}
In contrast to the backward differential flow based on Kuhn-Tucker's equation described in \cite{Zhu-2}, this paper introduces a method involving heat evolution. By applying heat conduct (Gaussian filtering) to the polynomial, we develop an algorithm similar to backward differential flow, where the evolution differential equation is explicitly formulated to aid in finding the minimizer. While the algorithm effectively reaches the global minimizer for quartic polynomials, it only partially succeeds with higher-degree polynomials. Such patterns were already observed in \cite{ABK}, demonstrated through Steklov regularization examples. However, we will present a theoretical criterion in Theorem \ref{suff-nece} that delineates the sufficient and necessary condition for these algorithms to achieve the global minimizer.

Our analysis is based on the Yuille-Poggio's fingerprints theory and their trajectory differential equation in the theory of computer vision \cite{Yuille-Poggio-1}\cite{Yuille-Poggio-2} which were built in 1980s. Before ending this introduction, we slightly sketch the motivation in our contributions.

 The elegant framework of multiscale Gaussian filter is equivalent to the model of heat conduct equation. Applying this theory, any even degree monic polynomials $p(x)$ will become convex by $p(x,t)=g_t(x)*p(x)$   for $t$ large enough\footnote{I definitely believe that this very simple fact should have been already established. But I have not gotten any references, limited to my scope of reading.}, where $g_t$ stands for Gaussian filter with variance $t$. Moreover,  for  quartic polynomial $p(x)$, the global
minimizer $x_{min}$  will continuously evolve along $t>0$ such that  it remains global minimizer $x^t_{min}$ of  $p(x,t)$ at each scale $t\ge 0$. Therefore,  reversely and continuously evolving from any global minimizer $x^t_{min}$ of $p(x,t)$ to $x_{min}$ of $p(x)$ is guaranteed. 

A natural question is, whether the global minimizer $x_{min}$ of a higher degree polynomial also evolves continuously to global minimizer of scaled version $p(x,t)$,  like the quartic polynomial case? Unfortunately this extremely expected  property doesn't hold in general for polynomials whose degree is more than $4$. We will  illustrate it by a counter-example on 6-degree polynomial. Furthermore, we give a condition which is both sufficient and necessary for the convergence to global minimizer.

 The multi-scale Gaussian filter and equivalent heat conduct equation is a standard content in the theory of PDEs, signal processing  and so on. In particular in the field of computer vision, it brought us many powerful theoretical tools since 1950s (cf. \cite{Iijima-1}\cite{Iikima-2}). Among them, the fingerprint theory proposed in 1980s (cf. \cite{Yuille-Poggio-1}\cite{Yuille-Poggio-2}) plays a kernel role for many years.

 In this paper, we apply the ideas of fingerprint from computer vision, and define three fingerprints of scaled polynomials $p(x,t)$  across scale $t$. The first fingerprint $\mathcal{FP}_1$ characterizes all the local extremals of $p(x,t)$ for each $t$, and the second one, $\mathcal{FP}_2$, characterizes the stationary points of $p(x,t)$ at each $t$, which indicate the domain of convexity of $p(x,t)$ during the time evolution. Furthermore, $\mathcal{FP}_3$ indicates the evolution of curves in $\mathcal{FP}_2$. All these  powerful fingerprints tools offer us insightful understandings to the evolution of both local and global extremals of the polynomials, from which we proposed a sufficient and necessary condition for attaining the global minimizer by the backward trajectory algorithm.

For the sake of simplicity, we list all the symbol and notations in this paper as Table \ref{list}.

\begin{table}[ht]
  \centering
  \begin{tabular}{c|p{5cm}|c}
  \hline
  Notations & Definition & Index\\
  \hline
  $\mathcal{Z}_{t,k}$ & real zeros set of $\frac{ \partial^k p(x,t)}{\partial x^k}$ & \eqref{eqn:zero_tk} \\
\hline
$\mu(x,t)$ &Steklov regularity of $p(x)$  & \eqref{eqn:steklov}\\
  \hline
  $p(x,t)$ & heat evolution of $p(x)$ & \eqref{eqn:heat} \\
  \hline
  $\mathcal{FP}_k(p)$ & $k$-th fingerprints, ($k=2$ is Yuille-Poggio's fingerprint) & \eqref{fp-1cd} \\
  \hline
  & Yuille-Poggio's fingerprint trajectory equation & \eqref{c-3}\\
  \hline
  $\mathbf{Flow}_{YP}(p)$ & Yuille-Poggio's flow & \eqref{fp-1cdabcx} \\
  \hline
  $\Omega(p)$ and $\Omega^c(p)$ & confinement zone and escape zone & Definition \ref{defn-1}\\
  \hline
$t_u$ & the merge time of finger prints for quartic polynomial & \eqref{single} \\
\hline
$t^*$ & the critical time after which the evoluted polynomial is convex&  \eqref{t-star}\\
\hline

    \end{tabular}
  \caption{List of notations and symbols}\label{list}
\end{table}

In conclusion, the primary contribution of this study is the identification of both necessary and sufficient conditions for the backward-differential-flow algorithm to achieve the global minimizer of a polynomial (Theorem \ref{suff-nece}). Our future work will explore the introduction of a seesaw PDE to secure the global minimizer, provided it resides within the defined {\it{attainable zone}}.

\section{Heat evolution and convexification of  polynomials}\label{Heat}
\subsection{Heat evolution of $p(x)$}\label{subsection-x}
Consider the heat conduct equation \eqref{a0} with initial condition \eqref{a0-ic}, the general solution of \eqref{a0} is
\begin{equation}\label{eqn:heat}
    p(x,t)=p(x)*g_t(x),
\end{equation}
in which $g_t(x)$ stands for the Gaussian filter
\begin{equation}\label{gauss}
    g_t(x)=\frac{1}{\sqrt{2\pi t}}e^{-\frac{x^2}{2t}},\ \ t\ge 0.
\end{equation}

In signal processing and computer vision, this time variable $t$ is also called {\it{scale}} (of Gaussian filtering) or {\it{artificial time}}. Notice that any differential operator  $\mathcal D$  is commutative with convolution operator $*$, i.e.,
\begin{equation}
{\mathcal D}(f*g)=f*{\mathcal D}g={\mathcal D}f*g.
\end{equation}

For polynomials $p(x,t)$, the heat equation \eqref{a0} can be enhanced to
\begin{equation}\label{a3}
\frac{\partial^k p}{\partial t^k}=\frac{1}{2^k}\frac{\partial^{2k} p}{\partial x^{2k}},\ \ \ \ (k=1,2,\cdots)
\end{equation}
by differentiating both sides of \eqref{a0} w.r.t $t$, since the smoothness is guaranteed. Then performing Taylor's expansion for $p(x,t)$ about  $t$ will yield
\begin{equation}\label{taylor-0}
p(x,t)=p(x,0)+t\cdot \frac{\partial p}{\partial t}\bigg |_{t=0}+\frac{t^2}{2}\frac{\partial^2 p}{\partial t^2}\bigg |_{t=0}+\cdots.
\end{equation}
If using the heat equation derived \eqref{a3}, we may rewrite \eqref{taylor-0} as
\begin{equation}\label{heat-1}
p(x,t)=p(x,0)+\frac{t}{2}\cdot \frac{\partial^2 p}{\partial x^2}\bigg |_{t=0}+\frac{t^2}{8}\frac{\partial^4 p}{\partial x^4}\bigg |_{t=0}+\cdots.
\end{equation}

The convexification  of  even degree polynomials  by heat evolution is characterized by following Theorem\footnote{Once again, I believe that this convexity result must be known in some literature.}.

\begin{theorem}\label{convex} For each even degree monic polynomial $p(x)$, there exists an specified $T^*=T^*(p)$ such that the heat convolution $p(x,t)$ is convex w.r.t $x$ at any $t>T^*$.
\end{theorem}




\begin{proof}
In what follows, the subscription $k$ in $P_k(x)$ and $Q_k(x)$ stands for the degree of polynomials. Let's consider even degree monic polynomial
\begin{equation}
P_{2m}(x)=x^{2m}+P_{2m-1}(x).
\end{equation}
Observing the expansion
\begin{equation}
P_{2m}(x)*g_t(x)=x^{2m}*g_t(x)+P_{2m-1}(x)*g_t(x),
\end{equation}
Clearly, 
we may  write
\begin{equation}
P_{2m}(x)*g_t(x)=P_{2m}(x)+\beta(x,t),
\end{equation}
in which
\begin{equation}
\beta(x,t)=(2m)!!\ t^m+ \sum_{k=1}^{m-1} t^{m-k}Q_{2k}(x).
\end{equation}
Using the heat evolution, we have
\begin{equation}
\frac{1}{2}\frac{\partial^2 p(x,t)}{\partial x^2}=\frac{\partial p(x,t)}{\partial t}=\frac{\partial \beta}{\partial t}.
\end{equation}
In our case,
\begin{equation}
\frac{\partial \beta}{\partial t}=m(2m)!!\ t^{m-1}+  \sum_{k=1}^{m-1} (m-k)t^{m-k-1}Q_{2k}(x).
\end{equation}

It is evident that all the primary terms of $Q_{2n}(x)$ originate from $x^{2m}(x)*g_t(x)-x^{2m}$ are necessarily positive. More specifically,
\begin{equation}\text{the leading coefficient of}\  Q_{2n}(x) =\binom{2m}{2n} (2m-2n)!! > 0.
\end{equation}
This suggests the existence of finite constants $K$, such that \begin{equation} Q_{2n}(x) > K > -\infty,\ \ (n=2,3,\ldots,2m-2).\end{equation}
Consequently,\begin{equation}\frac{\partial \beta}{\partial t} > m(2m)!!\ t^{m-1} + K( t^{m-2} + t^{m-3} + \ldots + 1).\end{equation} Thus, a value $T^* \in (0, +\infty)$ exists such that for all $t > T^*$, $\frac{\partial \beta}{\partial t} > 0$, indicating that convexity is maintained through heat evolution post $t > T^*$.

\end{proof}

\subsection{Comparison principle} The most important mechanism in heat evolution is the comparison principle, from which we understand that usually a local minimizer will merge to a local maximizer during the evolution, like the "annihilation" action between the pair of minimizer and maximizer.
\begin{theorem}[Comparison principle]\label{comparison-principle} Assume that $x^*$ be a critical point of $p(x,t^*)$, then for $t>t^*$, the heat evolution of the critical point satisfies
\begin{equation}
p( x^*(t),t)\ge p( x^*,t^*),\ {\rm{if}}\ x^*\ \rm{is\ local\ minimum};
\end{equation}
\begin{equation}
p( x^*(t),t)\le p( x^*,t^*),\ {\rm{if}}\ x^*\ \rm{is\ local\ maximum}.
\end{equation}
\end{theorem}

\begin{proof} Without loss of generality, we set $t^*=0$, due to that the heat operator $U^t: f(x)\mapsto g_t(x)*f(x)$ forms a semi-group (Lie group).
Let $x=x(t)$ be one of the integral curves of critical points of $p(x,t)$ w.r.t  $x$, then from
\begin{equation*}
\begin{split}
\frac{\,dp(x(t),t)}{\,dt}=&\frac{\partial p(x,t)}{\partial x}\dot{x}(t)+\frac{\partial p(x,t)}{\partial t}\\
 =&0+\frac{\partial p(x,t)}{\partial t}\\
 =&\frac{1}{2}\frac{\partial^2 p(x,t)}{\partial x^2},
\end{split}
\end{equation*}
thus we can get the required result. Notice that  the last equality comes from heat conduct equation.
\end{proof}


This comparative analysis indicates that local minimizers and local maximizers could potentially pair up and merge throughout the process. Ideally, an $n$ degree polynomial, where $n$ is even, should have $n-1$ critical points. Therefore, our expectation is that the global minimizer remains distinct from any local maximizers during the heat evolution. Nonetheless, there are scenarios where this may not hold, and we will examine these instances in depth.

\section{Global minimizer and scale space fingerprint}

\subsection{Fingerprints of scale space}
The scale space fingerprint  was introduced by A.L.~Yuille and T.A.~Poggio in 1980s (cf. \cite{Yuille-Poggio-1} \cite{Yuille-Poggio-2} etc.), which plays an important role in computer vision. In this scheme, the multi-scale version $p(x,t)=p(x)*g_t(x)$ of $p(x)$  comes from heat conduct equation, in which the variance $t\ge 0$ of Gaussian filter is also called artificial time.

Consider that all the polynomials in our situation are  of real coefficients, for the sake of simplicity, we need to generalize Yuille-Poggio's definition of fingerprints of multi-scale zero-crossings to more general case as below.

\begin{definition} Denote the set of real zeros of $k$-th derivative of polynomial $p(x,t)$ as
\begin{equation}\label{eqn:zero_tk}
 \mathcal Z_{t,k}(p):=\left \{ x_i(t)\in\mathbb R;\  \frac{\,\partial^k p(x_i(t),t)}{\,\partial x^k}=0,\ i=1,2,\cdots. \right \},
\end{equation}
and denote the  sets
\begin{equation}\label{fp-1c}
\begin{split}
\mathcal{FP}_k^+(p):=&\left \{ (x,t);\ \frac{\partial^k p(x,t)}{\partial x^k}>0,\ t\ge 0, \right \},\\
 \mathcal{FP}_k^-(p):=&\left \{ (x,t);\ \frac{\partial^k p(x,t)}{\partial x^k}<0,\ t\ge 0,\right \}.
\end{split}
\end{equation}
then  the $k$-th order fingerprints of polynomial $p(x)$ are defined as
\begin{equation}\label{fp-1ca}
\mathcal{FP}_k(p):=\overline{\mathcal{FP}_k^+(p)}\bigcap\overline{ \mathcal{FP}_k^-(p)}.
\end{equation}
\end{definition}

In above notations   $\overline {S}$ represents the topological closure of set $S$. 
In our case, this topological closure is very simple thus we may characterize $\mathcal{FP}_k$ by algebraic equations
 \begin{equation}\label{fp-1cd}
\mathcal{FP}_k(p)=\bigcup_{t:\ t\ge 0} \mathcal Z_{t,k}(p)=\left \{ (x,t);\ \frac{\partial^k p(x,t)}{\partial x^k}=0\right \},
\end{equation}
due to the sufficient smoothness of all polynomials.

\begin{remark} When $k=2$, the fingerprint $\mathcal{FP}_2$ of so-called \it{ zero-crossings}, as well as the equation of \it{zero-crossing contour},  are introduced by A.L. ~Yuille and T.A. ~Poggio   \cite{Yuille-Poggio-2}.  Here, we generalize their fingerprints from $\mathcal{FP}_2$ to more general $\mathcal{FP}_k$ ($k\ge 2$) in this paper. In other words, if we consider
$P(x)$ whose derivative is $P'(x)=p(x)$, then $\mathcal{FP}_1(p)=\mathcal{FP}_2(P)$. That is to say, our framework of $\mathcal{FP}_k$ is essentially a generalization of Yuille-Poggio's fingerprints  in the theory of computer vision.
\end{remark}

According to this notation, $\mathcal{FP}_1$ is the fingerprint of extremal values (critical points), and $\mathcal{FP}_2$ the zero-crossings (convexity)\footnote{Although there exists certain gap between the rigorous meaning and the definition here, we omit it in this paper.}, of polynomial $p(x)$, respectively. Essentially, as in the theory of computer vision, we can  get more information from
$\mathcal{FP}_2^+$ and $\mathcal{FP}_2^-$. In this paper, we generalize the classic concept $\mathcal{FP}_1$ and $\mathcal{FP}_2$ to general $\mathcal{FP}_k$, in particular, $\mathcal{FP}_3$ is included such that our main results can be represented on these three fingerprints.


We further consider the dynamics of the elements in $\mathcal {FP}_1$, i.e., the trajectories. Our main interest is to obtain the curves $x=x(t)$ which obey the equation
\begin{equation}\label{c-1}
\frac{\partial p(x(t),t)}{\partial x}=0,
\end{equation}
as well as initial conditions
\begin{equation}
 x(0)=x_i\in\mathcal{Z}_{0,1}(p),\ \ (i=1,2,\cdots)
\end{equation}
where $x_i$ ($i=1,2,\cdots$) are the critical points of $p(x)$.
To solve these curves, an ODE by varying $t$ as follows is introduced by A.L.~Yuille and T.A.~Poggio in \cite{Yuille-Poggio-2},
\begin{equation}\label{c-2}
0=\frac{d}{\,dt}\left ( \frac{\partial p(x(t),t)}{\partial x}   \right )=\frac{\partial^2 p(x,t)}{\partial x^2}\frac{\,dx(t)}{\,dt}+\frac{\partial^2 p(x,t)}{\partial x\partial t}.
\end{equation}
Therefore, we may characterize the fingerprint which contains all the maximums at different $t>0$ by rewriting \eqref{c-2} as
  \begin{equation}\label{c-3}
\frac{\,dx(t)}{\,dt}=-\frac{\frac{\partial^2 p(x,t)}{\partial x\partial t}}{\frac{\partial^2 p(x,t)}{\partial x^2}}=-\frac{\frac{\partial^3 p(x,t)}{\partial x^3} }{2\cdot\frac{\partial^2 p(x,t)}{\partial x^2}},
\end{equation}
as well as suitable initial conditions\footnote{If $p(x)$ is $n$-degree polynomial, there exists at most $n-1$ distinct initial conditions.}
\begin{equation}
x(0)\in \mathcal{Z}_{0,1}(p).
\end{equation}
In this paper, we call this ODE \eqref{c-3} the {\it{Yuille-Poggio equation}}, since it was first proposed in (3.3) of A.L.~Yuille and T.A.~Poggio's seminal work \cite{Yuille-Poggio-2}.

On the other hand, we also call this equation \eqref{c-3} the {\it trajectory equation}, since the reversely evolution algorithm will backward evolute along this curve, provided the initial value is given. Given any initial position, one may obtain a trajectory by this equation. In particular, when the initial condition is located at the critical points of $p(x)$, the trajectories form the fingerprint $\mathcal{FP}_1(p)$.

  We may further generalize $\mathcal{FP}_1(p)$ to  Yuille-Poggio's flow.

\begin{definition} For any $h\in\mathbb R$, the integral curve generated by Yuille-Poggio equation \eqref{c-3} associated with initial value $x(0)=h$ is called a Yuille-Poggio's curve. All these Yuille-Poggio's curves consist the set
\begin{equation}\label{fp-1cdabcx}
\mathbf{Flow}_{YP}(p):=\left \{ (x(t),t);\ \frac{\,dx(t)}{\,dt}=-\frac{\frac{\partial^3 p(x,t)}{\partial x^3}}{2\frac{\partial p^2(x,t)}{\partial x^2}},\ \   x(0)=h,\ \ \forall h\in \mathbb R   \right \},
\end{equation}
and we call it the Yuille-Poggio's flow generated by polynomial $p(x)$.
\end{definition}
Clearly, the fingerprint curve in the fingerprint $\mathcal{FP}_1(p)$ is a special Yuille-Poggio's curve whose  initial value $x(0)$  is restricted to $\mathcal{Z}_{0,1}(p)$, i.e., satisfies $p'(x(0))=0$.  Thus we have

\begin{theorem}\label{Y-P-theorem} The fingerprint $\mathcal{FP}_1$ can  be represented as
 \begin{equation}\label{fp-1cdabc}
\mathcal{FP}_1(p)=\left \{ (x,t);\ \frac{\,dx(t)}{\,dt}=-\frac{\frac{\partial^3 p(x,t)}{\partial x^3}}{2\frac{\partial p^2(x,t)}{\partial x^2}},\ \   x(0)\in \mathcal{Z}_{0,1}  \right \}.
\end{equation}
And
\begin{equation}
\mathcal{FP}_1(p)\subset \mathbf{Flow}_{YP}(p).
\end{equation}
\end{theorem}

Notice that the singularity occurs at which the denominator
vanishes, we have
\begin{lemma}\label{double-roots} For any polynomial $p(x)$ and its heat evolution $p(x,t)$,
if
\begin{equation}
(x',t')\in \mathcal{FP}_i\bigcap\mathcal{FP}_{i+1},
\end{equation}
then $x'$ must be a real double root of polynomial equation $\frac{\partial^i p(x,t)}{\partial x^i}=0$, and a real root of polynomial equation $\frac{\partial^{i+1} p(x,t)}{\partial x^{i+1}}=0$.
\end{lemma}

\begin{lemma} For $n$-th ($n$ is even) order polynomial $p(x)$, the set $\mathcal{FP}_2\bigcap \mathcal{FP}_3$ contains at most $\frac{n}{2}-1$ points $(x_i,t_i)$,
where $i=1,2,\cdots, \frac{n}{2}-1$.
\end{lemma}




\subsection{Confinement zone and escape zone}

In our following analysis, we will give basic framework of $\mathcal{FP}_2\bigcap\mathcal{FP}_3$,  we vary initial condition of trajectory ODE  to partition $\mathbb R$ into {\it Confinement Zone} and {\it Escape Zone}.

\begin{landscape}
    \begin{figure}
  \centering
\includegraphics[width=0.6\linewidth]{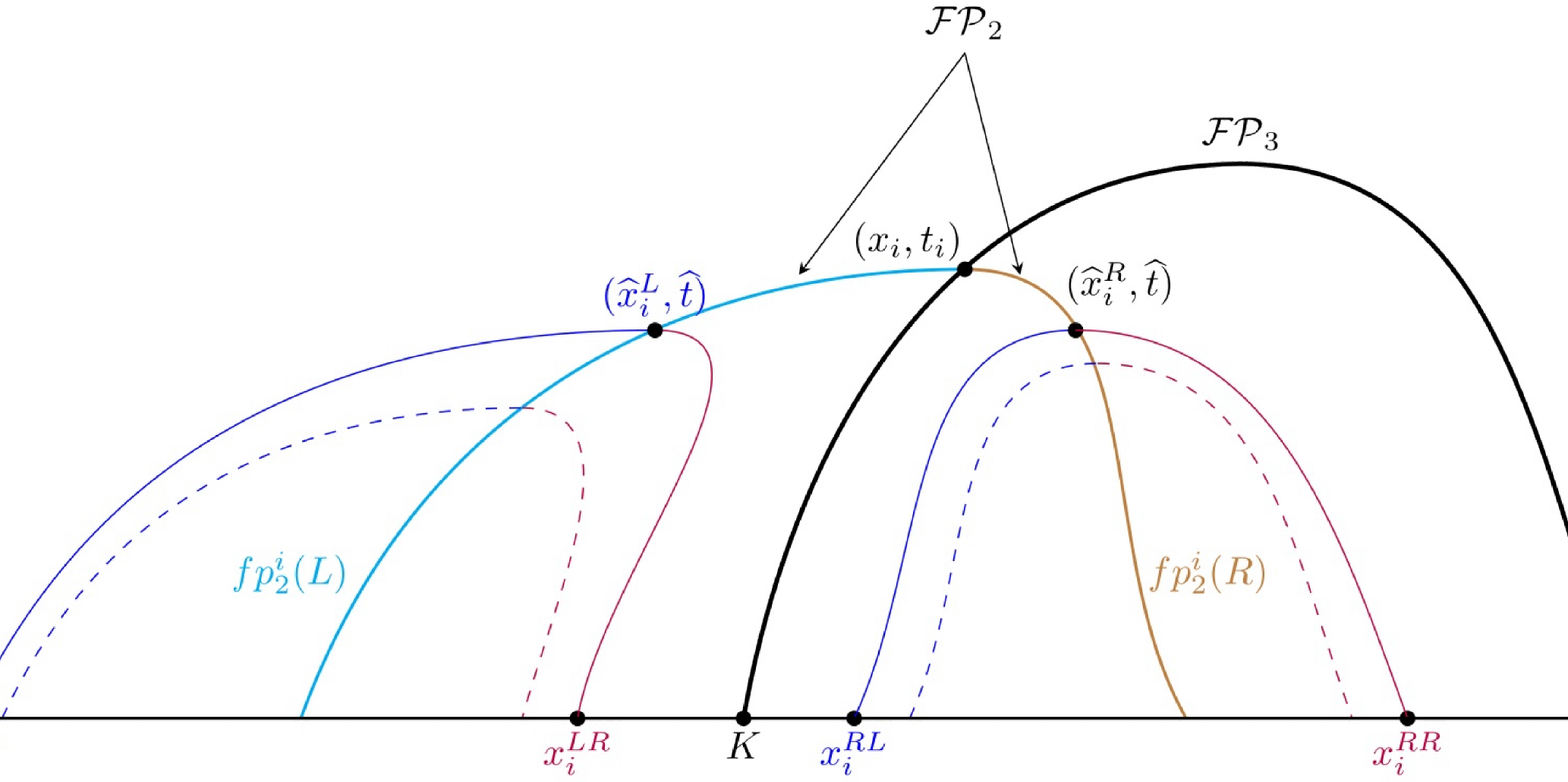}
  \caption{The illustration of Yuille-Poggio's flow as well as $\mathcal{FP}_2 $ and $\mathcal{FP}_3$.}
  \label{Yuille-Poggio's flow}
\end{figure}
\end{landscape}

\begin{definition}\label{defn-1} Let $c\in \mathbb R$, if the Yuille-Poggio's curve from $(c,0)$ will not intersect  with any Yuille-Poggio's curve from $(c',0)\ne (c,0)$, we call this $c$ is in {\it{Escape Zone}}. Otherwise, we say it is in the {\it{Confinement Zone}}, which is denoted by $\Omega$. Accordingly, the {\it{Escape\ Zone}} is denoted by $\Omega^c$.
\end{definition}

\begin{theorem}[Characterization of confinement zone and escape zone] The confinement zone $\Omega$ is
\begin{equation}\label{confinement-zone-1}
\Omega:=\bigcup_{i=1}^{\frac{n}{2}-1} [ X_i^{LL},X_i^{RR}].
\end{equation}
where
\begin{equation}\label{K-1}
X_i^{LL} =\lim_{L(\widehat x_i^L,\widehat t)\to (x_i,t_i)} x_i^{LL},
\end{equation}
\begin{equation}\label{K-4}
X_i^{RR}=\lim_{L(\widehat x_i^L,\widehat t)\to (x_i,t_i)} x_i^{RR},
\end{equation}
in which the limitation means the point $(\widehat x_i^L,\widehat t)$ (or $(\widehat x_i^R,\widehat t)$, resp.)  moves to the  destination $(x_i,t_i)$ along the local $\mathcal{FP}_2$ fingerprint curve $fp_2^i(L)$ (or $fp_2^i(R)$, resp.). Here
$(\widehat x_i^L,\widehat t)$ (or $(\widehat x_i^R,\widehat t)$, resp.) is the end of Yuille-Poggio curve connected to $(x_i^{LL},0)$ (or $(x_i^{RR},0)$, resp.).
 \end{theorem}

\begin{proof} We first prove that
\begin{equation}\label{confinement-zone-1a}
\Omega:=\bigcup_{i=1}^{\frac{n}{2}-1}\left ( [ X_i^{LL},X_i^{LR}]\bigcup [X_i^{RL},X_i^{RR}]\right ),
\end{equation}
in which we add two notations,
\begin{equation}\label{K-2}
X_i^{LR}=\lim_{L(\widehat x_i^L,\widehat t)\to (x_i,t_i)} x_i^{LR}=K,
\end{equation}
\begin{equation}\label{K-3}
X_i^{RL}=\lim_{L(\widehat x_i^L,\widehat t)\to (x_i,t_i)} x_i^{RL}=K.
\end{equation}
 Here $K$ stands for the  intersection point $(K,0)$ between the curve in $\mathcal{FP}_3(p)$ and straight line $t=0$.

Let's show that the right hand side of \eqref{confinement-zone-1a} is well defined. As illustrated at Fig.\ref{Yuille-Poggio's flow}, connecting to each $(x_i,t_i)\in {\mathcal{FP}_2\bigcap\mathcal{FP}_3}$, there exists a pair of curves in $\mathcal{FP}_2$, corresponding to $(x_i+0,t_i-0)$ and $(x_i-0,t_i-0)$,  denoted by  $fp_2^i(R)$ and $fp_2^i(L)$ respectively.

  For any point $(\widehat x_i^R,\widehat t)\in fp_2^i(R)$, when $(\widehat x_i^R,\widehat t)\ne (x_i,t_i)$, there are a pair of trajectories satisfying \eqref{c-3} which contains  the point $(\widehat x_i^R,\widehat t)$. We may denote their  ends at $t=0$ as $(x_i^{RL},0)$ and $(x_i^{RR},0)$ respectively. Here we assume that $x_i^{RL}\le  x_i^{RR}$.

  Similarly, for any point  $(\widehat x_i^L,\widehat t)\in fp_2^i(L)$, when $(\widehat x_i^L,\widehat t)\ne (x_i,t_i)$, there are a pair of trajectories satisfying \eqref{c-3} which contains  the point $(\widehat x_i^L,\widehat t)$. We denote their ends at $t=0$ as $(x_i^{LL},0)$ and $(x_i^{LR},0)$ respectively. Here we assume that $x_i^{LL}\le x_i^{LR}$.

     Now we may write that
  \begin{equation}x_i^{LL}\le x_i^{LR}<K< x_i^{RL}\le x_i^{RR}\end{equation}
  due to that the Yuille-Poggio curve can not intersect with $\mathcal{FP}_3(p)$ otherwise it will bring singularities, according to the denominator of the right hand side of Yuille-Poggio equation \eqref{c-3}.

  Furthermore, the limitation process in \eqref{K-1} etc. remains  monotonicity. That is, moving from $({\widehat x_i}^L,\widehat t)$ to $({\widehat x_i}^{L+},{\widehat t}^+)$  and finally to  $(x_i,t_i)$, we may observe that
  \begin{equation} -\infty<\widehat x_i^{LL+}<   \widehat x_i^{LL}<\widehat x_i^{LR} <\widehat x_i^{LR+}<K.\end{equation}
  This implies that all the limitation \eqref{K-1} and so on are well-defined.  Finally, we note that  $\forall h\in (K-\epsilon,K)$, there must exist a
  Yuille-Poggio curve starts from $(h,0)$, and for $\forall \epsilon>0$, for any point
  $  (x',t')\in fp_2^i(L)$, that satisfy $ t'<t, x'>x_i$ and $\|(x',t')-(x_i,t_i)\|_2<\epsilon$, there must exist a Yuille-Poggio curve pass the point $(x',t')$.
  That is to say, the Yuille-Poggio curves near the $\mathcal{FP}_3$ curve connecting $(K,0)$ and $(x_i,t_i)$, are dense. Here for the sake of simplicity, we omit the  topology and differential dynamics description.

  Now  the set $\Omega$ in \eqref{confinement-zone-1} is well defined. We observe that any Yuille-Poggio curve starting from $(h,0)$ for $h\in\Omega$ occurs if and only if there exists another Yuille-Poggio curve, starting from $(h',0)$ for some $h'\in \Omega$. In particular, these two curves meet at $fp_2^i(L)$ or $fp_2^i(R)$. Thus the current $\Omega$ in \eqref{confinement-zone-1a} is agreed with that in Definition \ref{defn-1}.

  \end{proof}

Now we state the main contribution of this paper:

\begin{theorem}\label{suff-nece} Let $p(x)$ be  polynomial of even order with a positive leading coefficient. Suppose $x^{t_0}_{\rm{min}}$ is the global minimizer for the convex polynomial (appropriately scaled version) $p(x,t) = p(x)*g_t(x)$ at $t = t_0$, and the trajectory's endpoint according to Yuille-Poggio's trajectory equation \eqref{c-3} at $t = 0$ is $x^*$. Then the global minimizer $x^*$ of $p(x)$ can be inversely involved from the global minimizer of its convexification version $p(x,t_0)$, {\bf{if and only if}} $x^*$ is in the \it{Escape\ Zone} $\Omega^c$.
\end{theorem}
\begin{proof} If optimal minimizer $x^*\in\Omega$,  then the maximum of $t$ associated with all the Yuille-Poggio curves is bounded, thus all these Yuille-Poggio curves can not connect to the point in $\mathbb R_x\times \mathbb R_t$ with  large $t>0$. Henceforth, the global minimizer $x^*\in \Omega^c$. 
\end{proof}

\begin{remark}
While it is anticipated that $X_i^{LL},X_i^{LR},X_i^{RL},X_i^{RR}$ would be explicitly defined, their algebraic expressions are unattainable because for polynomials of degree $6$ or higher, the corresponding curve in $\mathcal{FP}_1$ necessitates dealing with algebraic equations of at least degree $5$. Consequently, we propose employing numerical techniques to determine these values.
\end{remark}

\begin{remark}
The analytical approach presented for convexification through the heat conduct equation, specifically Gaussian filtering, remains applicable to Steklov regularization scenarios.
\end{remark}

\section{Case study of quartic polynomials}

In what follows, we will give explicit representation for the fingerprints of quartic polynomials, and explain their geometric properties. Although partial results are the similar as \cite{ABK} for quartic polunomials, we found a distinct algorithm to attain the global minimizer of quartical polynomials without any convexification treatment, neither Steklov's regularization nor heat evolution.

\subsection{The structure of fingerprints}

For the quartic polynomial 
\begin{equation}\label{eq:no1} p(x)=x^4+ax^3+bx^2+cx+d,\end{equation}
we see its heat evolution
\begin{equation}\label{b-0}
\begin{split}
    p(x,t)&=p(x)+ (6x^2+3ax+b)\cdot t+3t^2\\
    &=x^4+ax^3+(b+6t)x^2+(c+3at)x+(d+bt+3t^2).
\end{split}
\end{equation}

Continue to differentiate  both sides of \eqref{b-0} w.r.t $x$,  the information of $\frac{\partial p(x,t)}{\partial x}$
across time $t$ may be represented as
 \begin{equation}\label{b-1}
    \frac{\partial p(x,t)}{\partial x}=4x^3+3ax^2+(2b+12t)x+(c+3at) = 0.\\
\end{equation}
Similarily, we have
 \begin{equation}\label{b-2}
    \frac{\partial^2 p(x,t)}{\partial x^2}=\frac{\partial^2 p(x)}{\partial x^2}+ 12 t=(12x^2+6ax+2b)+12t=0,
\end{equation}
and
 \begin{equation}\label{b-2-1}
    \frac{\partial^3 p(x,t)}{\partial x^3}=24x+6a=0.
\end{equation}
These three equations form the characterization of fingerprints $\mathcal{FP}_1$, $\mathcal{FP}_2$ and $\mathcal{FP}_3$, respectively.

\subsubsection{The structure of fingerprint $\mathcal{FP}_1$}

Based on \eqref{b-1}, the fingerprint $\mathcal {FP}_1$ is characterized by following time-varying cubic equation
\begin{equation}\label{eqn:fp-1c}
x^3+\frac{3a}{4} x^2 + \frac{b+6t}{2} x+\frac{c+3at}{4}=0,
\end{equation}
Now, if $x_i$ is a real root of \eqref{eqn:fp-1c} at $t=0$, then it leads to the trajectory described by the differential equation \eqref{c-3}. For more details, we have

\begin{lemma} For quartic polynomial $p(x)$, the local extremal values points $x^t_i\ (i=1,2,3)$ of $p(x,t)$ w.r.t $x$ at scale $t$ satisfy the trajectory differential equation
\begin{equation}\label{eqn:dyn-1d}
\frac{\,d{ x}(t)}{\,dt}=-\frac{12{x}+3a}{12{x}^2+6ax+2b+12t},
\end{equation}
with following (at most three) initial conditions,
\begin{equation}
x_i(0)= x_i,\ \ \ (i=1,2, 3).
\end{equation}
Here $ x_i$ is  the local extremal of $p(x)$.
\end{lemma}
\begin{proof} 
Inserting \eqref{b-2} and \eqref{b-2-1} into \eqref{c-3} will lead to required results.

\end{proof}

The equation \eqref{eqn:fp-1c} possesses (at most) three real roots at $t=0$, corresponds to (at most) three trajectories, which form the Fingerprint $\mathcal{FP}_1$.
However, on the viewpoint of differential algebra (see, e.g. \cite{Kap}), actually the solution of differential equation \eqref{eqn:dyn-1d} is real algebraic curve,
i.e., a polynomial $F(x,t)$ about $x(t)$ and $t$ which satisfy $F(x,t)=0$. In our case, the polynomial equation \eqref{eqn:fp-1c} describes this algebraic curve, thus we may immediately apply the algebraic representation of $\mathcal{FP}_1$:
\begin{equation}\label{eqn:algebraic-curve}
\mathcal{FP}_1=\left \{(x,t);\ t=-\frac{4x^3+3ax^2+2bx+c}{12x+3a},\ \ x\ne -\frac{a}{4},\ {\rm{and}}\ t\ge 0\right\}.
\end{equation}

 According to Subsection \ref{discriminant}, to get the information of the  roots of \eqref{eqn:fp-1c}, we need  its  discriminant,
\begin{equation}
\Delta(t)=\left ( \frac{a^3-4ab+8c}{64}\right )^2 + \left ( \frac{-3a^2+8b}{48}+t\right )^3,
\end{equation}
which will be explained in details in \eqref{eq:d-2}.

\begin{lemma}\label{lem:discriminant} The discriminant $\Delta(t)$ of equation \eqref{eqn:fp-1c} is monotonically increasing to infinity. Its unique zero is
\begin{equation}\label{single}
t_u=\frac{a^2}{16}-\frac{b}{6}-\frac{1}{16}(a^3-4ab+8c)^{\frac{2}{3}}.
\end{equation}
\end{lemma}

\begin{proof}
Using \eqref{eq:d-1}, we write
\begin{equation}\label{eqn:cubic-f-g}
\begin{split}
f(t)&=\frac{b}{2}-\frac{3a^2}{16}+3t,\\
g(t)&=\frac{a^3}{32}-\frac{ab}{8}+\frac{c}{4}.\\
\end{split}
\end{equation}
Now  the time-varying discriminant
\begin{equation}\label{eq:d-2}
\begin{split}
\Delta(t)=&\frac{[g(t)]^2}{4}+\frac{[f(t)]^3}{27},\\
=&\left (\frac{a^3-4ab+8c}{64} \right )^2+\left (\frac{-3a^2+8b}{48}+t \right )^3,
\end{split}
\end{equation}
which means that  $\Delta(t)$ increases monotonically w.r.t. $t$. Immediately, \eqref{eq:d-2} leads to \eqref{single}.
\end{proof}

\begin{theorem}[The "$1+2$" structure of $\mathcal{FP}_1$]\label{quartic-fingerprints}
Let $t_u$, defined in \eqref{single}, be the zero of discriminant $\Delta(t)$. If $t_u< 0$, then $\mathcal{FP}_1$ contains only one trajectory $x(t)$ described by equation \eqref{eqn:dyn-1d}, which evolutes as $t\to +\infty$. If $t_u \ge 0$, then during $t\in [0,t_u]$ the Fingerprint $\mathcal{FP}_1$ contains three distinct trajectories described by \eqref{eqn:dyn-1d}, one of which continues to evolute to $+\infty$ during $t>t_u$, and the other two  trajectories will start from $t=0$ but merge (stop) when $t=t_u$ at the point $(x(t_u),t_u)$. Here,
 \begin{equation}\label{singularity-1}
x(t_u)=\left (\frac{a^3-4ab+8c}{64} \right )^{1/3}-\frac{a}{4}.
\end{equation}
\end{theorem}
\begin{proof} According to Lemma \ref{lem:discriminant}, we know that if $t_u<0$, then $\Delta(t)>0$ for all $t\ge 0$, which means the equation \eqref{eqn:fp-1c} has only one root at each $t\ge 0$. When $t_u\ge 0$, then $\Delta(t)<0\ (=0,\  >0,\ respectively)$ while $t\in [0,t_u)$ ($t=t_u$,\ $t>t_u$, respectively), and the equation \eqref{eqn:fp-1c} has three distinct real roots (one real and a pair of double real roots, or one real root, respectively).
 In particular, we consider the critical case
$t=t_u$ at which $\Delta(t)=0$.
If case is this, the equation \eqref{eqn:fp-1c} at $t=t_u$ possesses one real root and a real double root.
It follows from \eqref{eqn:cubic-3} that the real double root is
\begin{equation}
  x(t_u)=\left (\frac{g(t)}{2} \right )^{\frac{1}{3}}-\frac{1}{3}\cdot\frac{3a}{4}.
\end{equation}
Substituting  \eqref{eqn:cubic-f-g} into this formula will produces \eqref{singularity-1}.
\end{proof}

\subsubsection{The structure of $\mathcal{FP}_2$ and $\mathcal{FP}_3$}

The structure of $\mathcal{FP}_3$ is very simple for quartic polynomial, since from \eqref{b-2-1} we may write
\begin{equation}
\mathcal{FP}_3=\left \{(x,t);\ x=-\frac{a}{4},\ t\ge 0\right \}.
\end{equation}

To analyze the structure of $\mathcal{FP}_2$, we have a Lemma as below.
\begin{lemma} Denote
\begin{equation}\label{t-star}
 t^*=\frac{a^2}{16}-\frac{b}{6},
\end{equation}
then the polynomial $p(x,t)$ defined in \eqref{b-0} is convex about $x$ at each   $t>\max(t^*,0)$. Furthermore, this $t^*$ can not be improved.
\end{lemma}
\begin{proof}
Consider the lower bound of  \eqref{b-2} at $t=0$,
\begin{equation}\label{b-3}
\begin{split}
\frac{\partial^2 p(x)}{\partial x^2}&=12x^2+6ax+2b\\
&=12\left (x+\frac{a}{4}\right )^2-\frac{3a^2}{4}+2b\\
&\ge -\frac{3a^2}{4}+2b=-12t^*.
\end{split}
\end{equation}
Combining it with \eqref{b-2}, we would have
\begin{equation}\label{b-4}
\frac{\partial^2 p(x,t)}{\partial x^2}\ge -\frac{3a^2}{4}+2b+12t=12(t-t^*),
\end{equation}
which implies the required results.

On the other hand, at any fixed $t'<t^*$, notice that  at $x=-\frac{a}{4}$, we have
\begin{equation}\label{b-3hh}
\begin{split}
\frac{\partial^2 p(x,t')}{\partial x^2}&=12x^2+6ax+2b+12t'\\
&=12\left (x+\frac{a}{4}\right )^2-12(t^*-t')\\
& = -12(t^*-t')<0,
\end{split}
\end{equation}
which is not convex at this $x=-\frac{a}{4}$, such that  $t^*$ is the optimal, and can not be improved.
\end{proof}

\begin{theorem}[The structure of $\mathcal{FP}_2$] For the structure of $\mathcal{FP}_2$ of quartic polynomial $p(x)$,
 \begin{enumerate}
                  \item when $t^*<0$, the fringerprint $\mathcal {FP}_2$ is empty.
                                                              \item when $t^*=0$, the $$\mathcal{FP}_2=\left \{(x,t)=(-\frac{a}{4},0)\right \}$$ has only single element;
                  \item when $t^*>0$,
the fingerprint $\mathcal{FP}_2$  consists of two curves: the left one is
\begin{equation}
x_L(t)=-\frac{a}{4}- \sqrt{t^*-t},\ \ \ (t^*\ge t\ge 0),
\end{equation}
and the right one is
\begin{equation}
x_R(t)=-\frac{a}{4}+ \sqrt{t^*-t},\ \ \ (t^*\ge t\ge 0).
\end{equation}
Specifically, these two curves must meets up at $t=t^*$, i.e., at the point
\begin{equation}
(x_L(t^*),t^*)=(x_R(t^*),t^*)=\left  (-\frac{a}{4},t^*\right ).
\end{equation}
\end{enumerate}
\end{theorem}

\begin{proof}    (1) comes from the fact  that for every $t\ge 0$, all the $\frac{\,\partial^2 p}{\,\partial x^2}>0$. That is, $\mathcal{FP}_2^+=\{(x,t);\ x\in\mathbb R,\ t\in [0,+\infty)\}$, but $\mathcal{FP}_2^-=\emptyset$. (2) is an immediate result, and (3) is from the quadratic equation  \eqref{b-4}.
 \end{proof}

\subsubsection{The intersection between fingerprints}

According to  Lemma \ref{double-roots}, we may summarize the intersection of fingerprints.

\begin{theorem} For the monic  quartic polynomials $p(x)=x^4+ax^3+bx^2+cx$, the  two intersection sets
\begin{equation}
\mathcal{FP}_2\bigcap\mathcal{FP}_3=\left \{(-\frac{a}{4},t^*)\right \},
\end{equation}
and
\begin{equation}
\mathcal{FP}_1\bigcap\mathcal{FP}_2 =\{(x(t_u),t_u)\},
\end{equation}
in which $t^*$ is defined in \eqref{t-star}, $t_u$ and $x(t_u)$ are defined in \eqref{single} and \eqref{singularity-1} respectively.
\end{theorem}

\begin{remark}[Three phase of time evolution] In general settings, the evolution of polynomial $p(x)$ can be categorized into three phases according to $0\le t_u\le t^*$. At first phase, $t$ evolutes from $0$ to $t_u$, and $\mathcal{FP}_1$ contains three distinct trajectories. Two of them will merge at $t=t_u$.

Then at the second phase, $t_u<t<t^*$, the $\mathcal{FP}_1$ contains only one trajectory, but $p(x,t)$ is not convex.

Finally, at the third phase, $t>t^*$, the $\mathcal{FP}_1$ contains only one trajectory, and $p(x,t)$ is  convex.

\end{remark}

\subsection{Confinement zone}

Now we compute the confinement zone. 
\begin{theorem} The confinement zone of quartic polynomial $p(x)$ is the interval
\begin{equation}
\Omega=\left [-\frac{a}{4}-\sqrt{3t^*},-\frac{a}{4}+\sqrt{3t^*} \right ],
\end{equation}
where $t^*$ is defined in \eqref{t-star}.
\end{theorem}
\begin{proof} 
We rewrite \eqref{single} as
\begin{equation}\label{single-star}
t_u(c)=t^*-\frac{1}{16}(a^3-4ab+8c)^{\frac{2}{3}}.
\end{equation}
Clearly, We should vary $c$ such that $t_u(c)=t^*$, i.e.,
\begin{equation}
a^3-4ab+8c=0\implies c=\frac{ab}{2}-\frac{a^3}{8}.
\end{equation}
Substituting this $c$ into the trajectory algebraic curve equation \eqref{eqn:algebraic-curve} and setting $t=0$, we get the equation
\begin{equation}
4x^3+3ax^2+2bx+ \frac{ab}{2}-\frac{a^3}{8} =0.
\end{equation}
The three roots of this equation are
\begin{equation}
x_1=-\frac{a}{4},\ x_{2,3}=-\frac{a}{4}\pm\sqrt{3t^*},
\end{equation}
which produces two pair of trajectories started from $t^*$ but reversely evolutes to $t=0$, whose four destinations form the confinement zone
\begin{equation}
\left [-\frac{a}{4}-\sqrt{3t^*},-\frac{a}{4} \right ]\bigcup \left [-\frac{a}{4},-\frac{a}{4}+\sqrt{3t^*} \right ]
=\left [-\frac{a}{4}-\sqrt{3t^*},-\frac{a}{4}+\sqrt{3t^*} \right ]
\end{equation}
\end{proof}


\subsection{Differential equation of critical points across scale}


Now we give the representation of $t^*$ and $t_u$ in terms of roots of fingerprint equation $\frac{\partial p(x,t)}{\partial x}=0$.

\begin{lemma}\label{13} Let $x_1,x_2,x_3$ be the critical points of \eqref{a1}, $t^*$ is defined in \eqref{t-star}, and $t_u$ defined in \eqref{single}, then we have
\begin{equation}
t^*=\left (\frac{x_1+x_2+x_3}{3}\right )^2-\frac{x_1x_2+x_2x_3+x_3x_1}{3},
\end{equation}
and
\begin{equation}\label{t-start-t-u}
t_u=t^*-\left [\frac{32}{27}(2x_1-x_2-x_3)(2x_2-x_3-x_1)(2x_3-x_1-x_2) \right ]^{\frac{2}{3}}.
\end{equation}
\end{lemma}
\begin{proof}
$x_1,x_2,x_3$ satisfy $p'(x)=0$, i.e.,
\begin{equation}\label{basic-2}
x^3+\frac{3a}{4}x^2+\frac{b}{2}x+\frac{c}{4}=0.
\end{equation}
thus
\begin{equation}\label{abc-by-x123}
\begin{split}
a=&-\frac{4}{3}(x_1+x_2+x_3),\\
 b=&2(x_1x_2+x_2x_3+x_3x_1),\\
  c=&-4x_1x_2x_3.\\
  \end{split}
\end{equation}
Substituting this into \eqref{single} and \eqref{t-star} will produce all these results. 

\end{proof}

\begin{theorem}\label{singularity-xt}
The singularity of the equation \eqref{eqn:dyn-1d} occurs only at
\begin{equation}\label{singularity}
x^{t_u}=x(t_u)=\left (\frac{a^3-4ab+8c}{64} \right )^{1/3}-\frac{a}{4},\ \ t=t_u.
\end{equation}
\end{theorem}
\begin{proof}
The singularity occurs at differential equation \eqref{eqn:dyn-1d}, which describes the $\mathcal{FP}_1$, so it must satisfy \eqref{eqn:fp-1c}. Meanwhile, the denominator of the r.h.s. of \eqref{eqn:dyn-1d} is actually the fingerprint of $\mathcal{FP}_2$, which should satisfy \eqref{b-4}. Thus we may combine these two algebraic equations to solve $(x,t)$.

Multiplying both sides of \eqref{b-4} by $x+\frac{a}{4}$, and subtracted it from \eqref{eqn:fp-1c} will produce
\begin{equation}\label{happy}
t=\frac{(3a^2-8b)x-(6c-ab)}{48x+12a},\ \ \ x\ne -\frac{a}{4}.
\end{equation}
Substituting it into \eqref{b-4} will yield a cubic equation about $x$,
\begin{equation}
48x^3+36ax^2+9a^2x+(3ab-6c)=0.
\end{equation}
This cubic equation has only one real solution  \eqref{singularity}. Substituting this $x$ into
\eqref{happy} will show that $t=t_u$.
\end{proof}

When $-\frac{a}{4}$ does not serve as a critical point, the point $(x^{t_u},t_u)$ is found exclusively on two $\mathcal{FP}_1$ integral curves from \eqref{eqn:dyn-1d}, which originate from a local minimum and a local maximum, respectively. One of these curves is associated with $\dot x(t_u)=+\infty$, while the other corresponds to $\dot x(t_u)=-\infty$. An interesting observation is that the integral curve that begins at a global minimum does not encounter this point $(x^{t_u},t_u)$, a detail that will be elaborated upon in Section \ref{Good}.

\begin{figure}[ht]
\centering

\begin{subfigure}[c]{0.45\linewidth}
\centering
\includegraphics[width=1.3\linewidth]{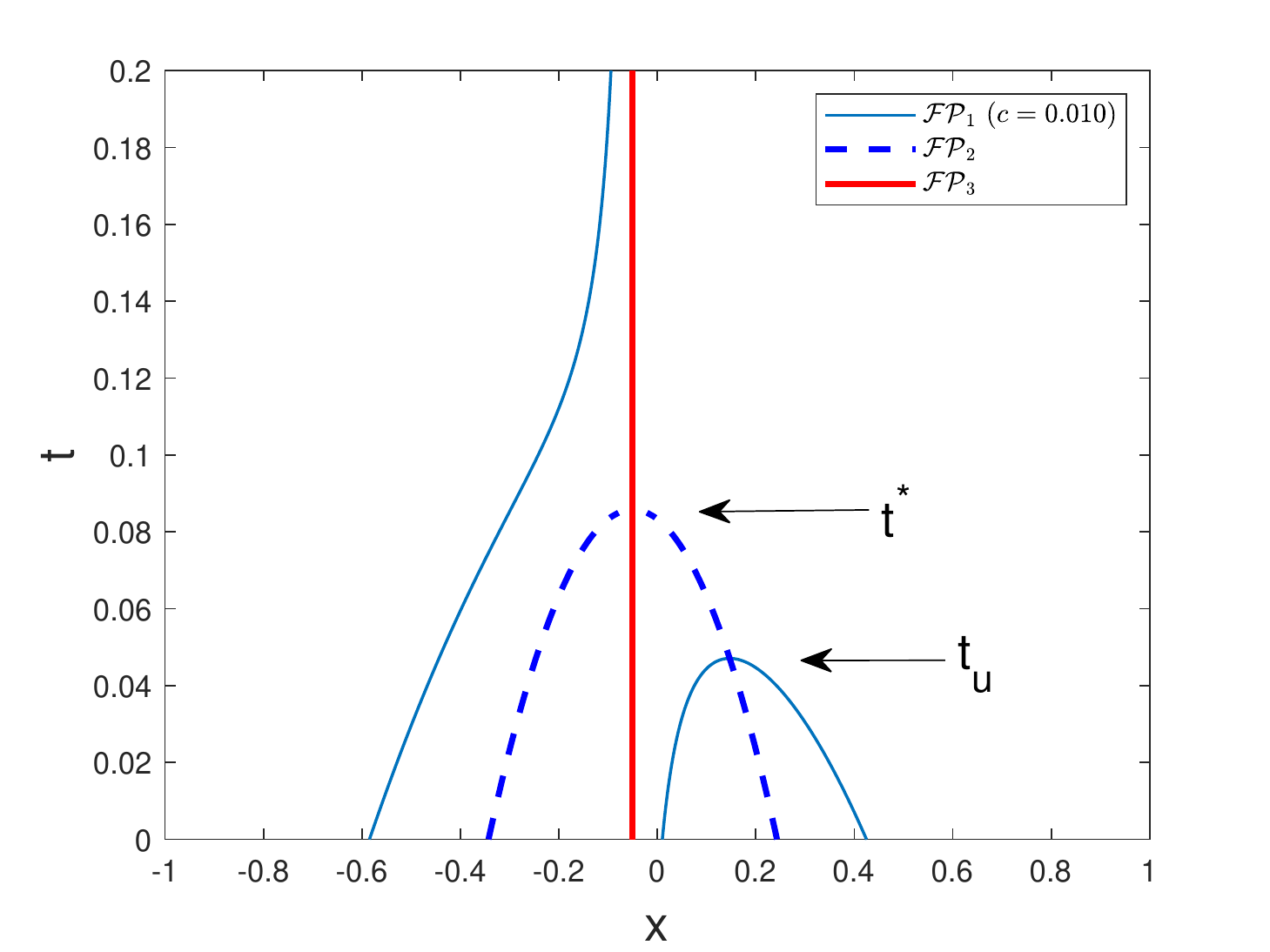}
  \caption{}
    \label{fig:sub---a}
\end{subfigure}
\hfill 
\begin{subfigure}[c]{0.45\linewidth}
    \centering
    \includegraphics[width=1.3\linewidth]{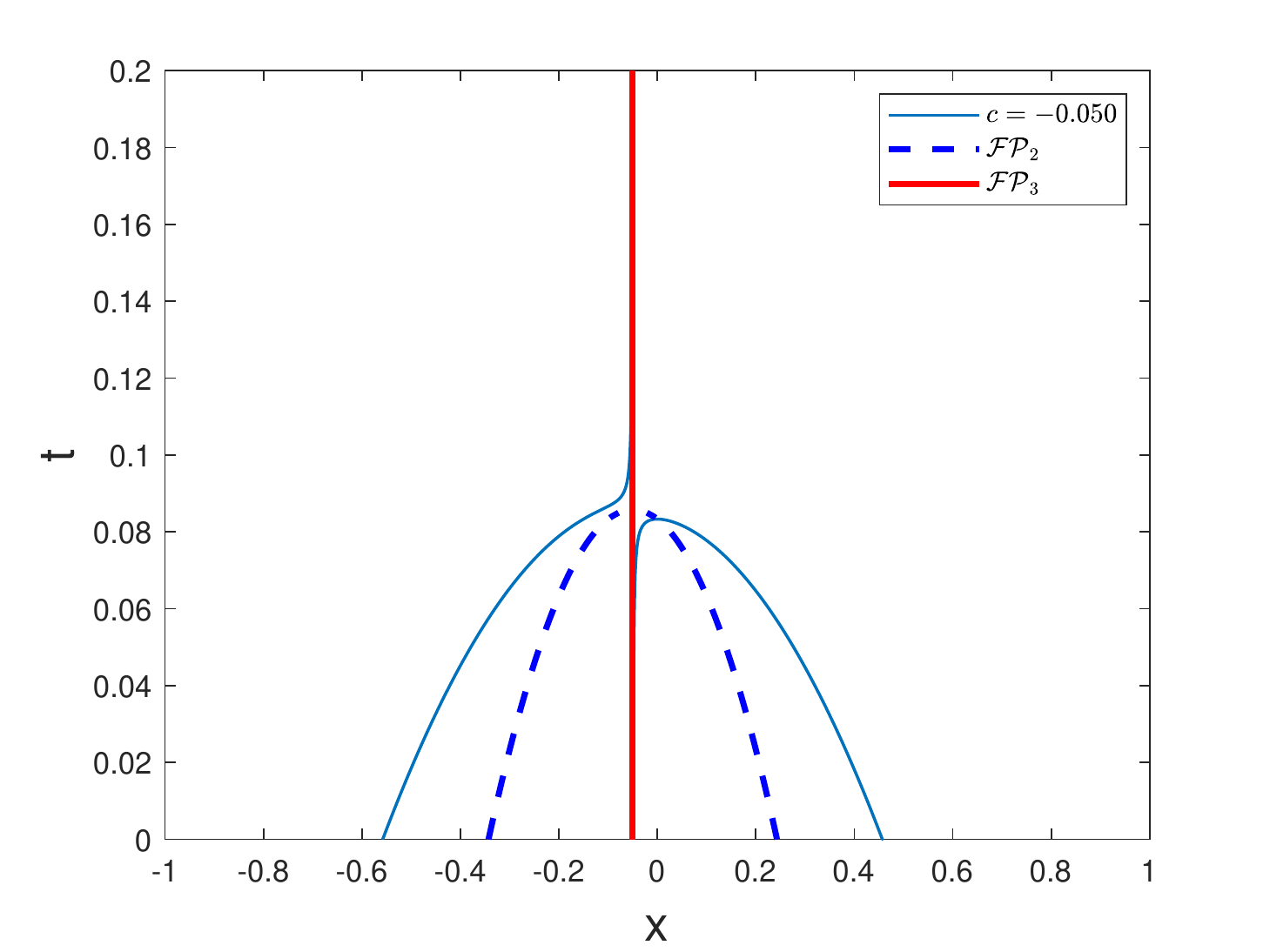}
   \caption{}
    \label{fig:sub---b}
\end{subfigure}
\begin{subfigure}[c]{0.45\linewidth}
\centering
\includegraphics[width=1.3\linewidth]{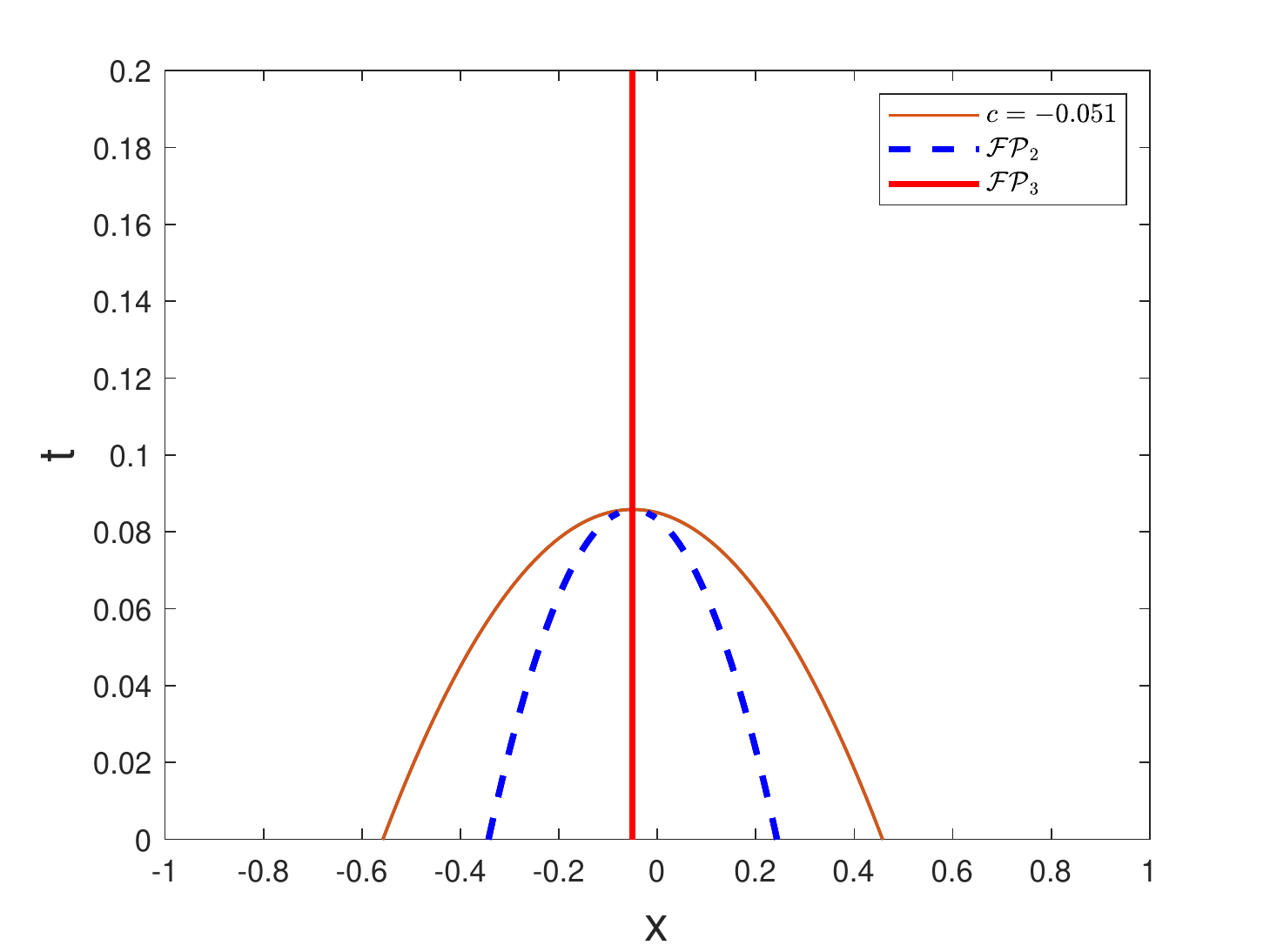}
\caption{}
    \label{fig:sub:b}
\end{subfigure}
\hfill
\begin{subfigure}[c]{0.45\linewidth}
\centering
\includegraphics[width=1.3\linewidth]{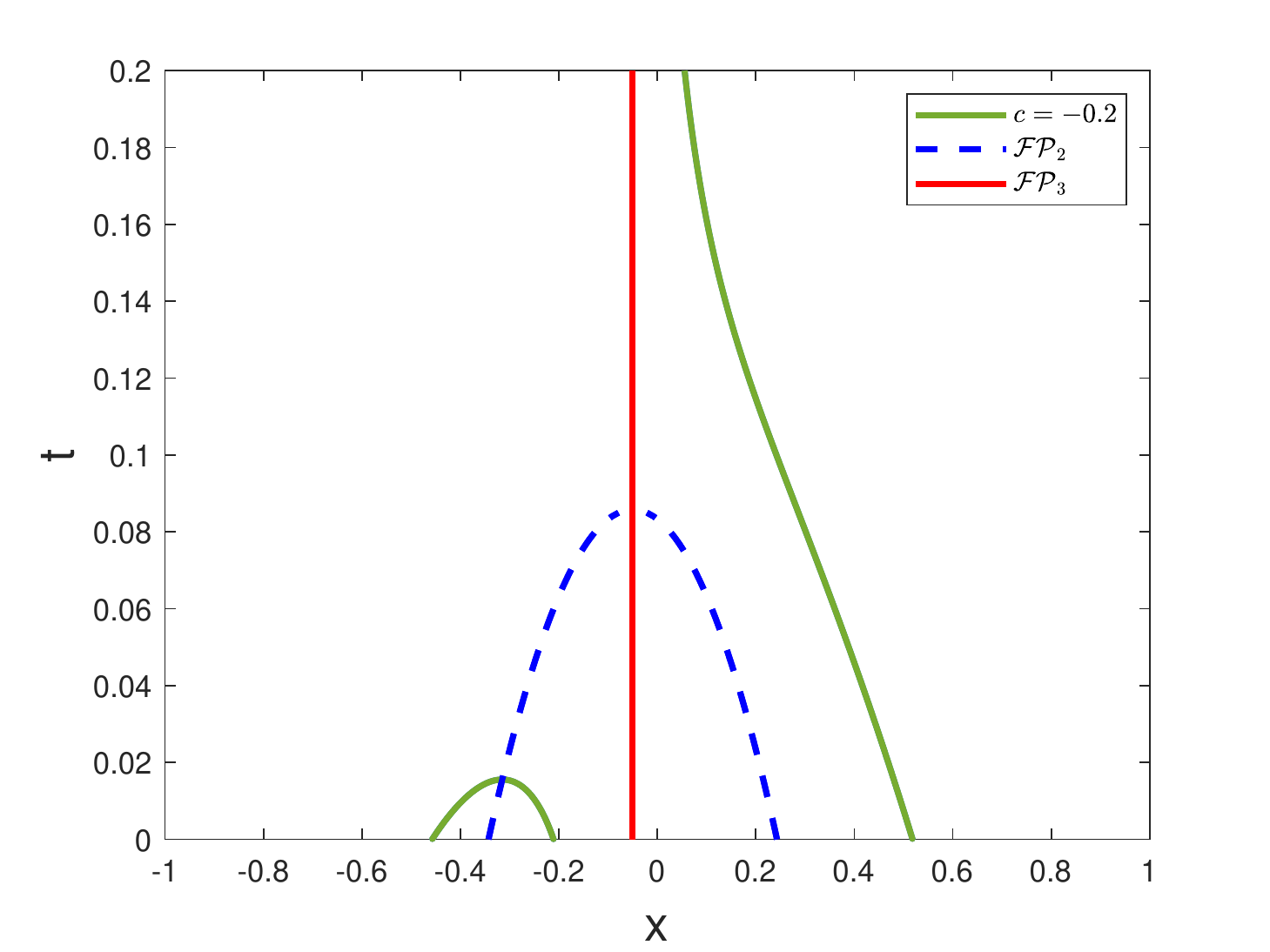}
\caption{}
\end{subfigure}

\begin{subfigure}[c]{0.45\linewidth}
\centering
\includegraphics[width=1.3\linewidth]{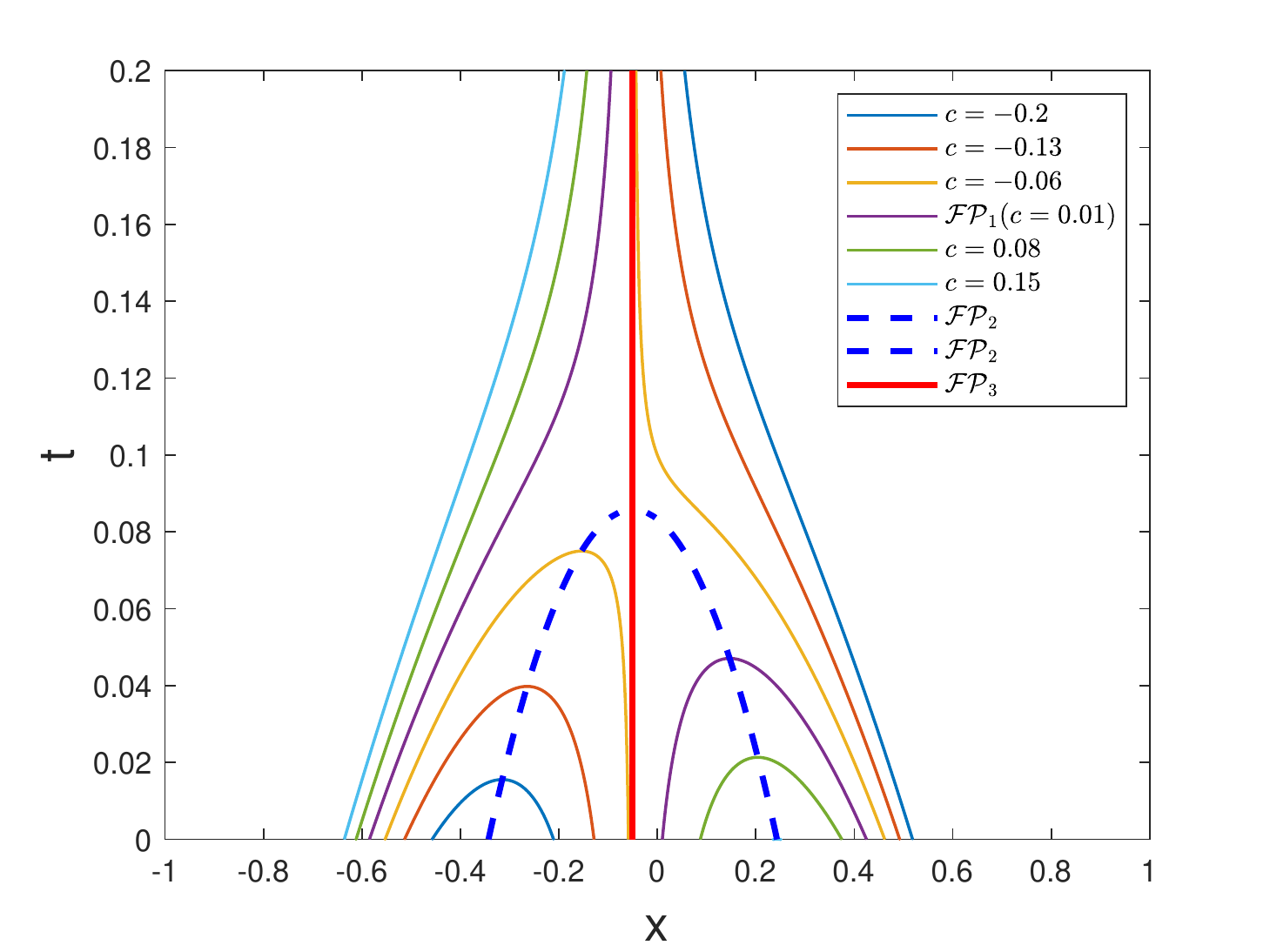}
\caption{}
\end{subfigure}
\hfill
\begin{subfigure}[c]{0.45\linewidth}
\centering
\includegraphics[width=1.4\linewidth]{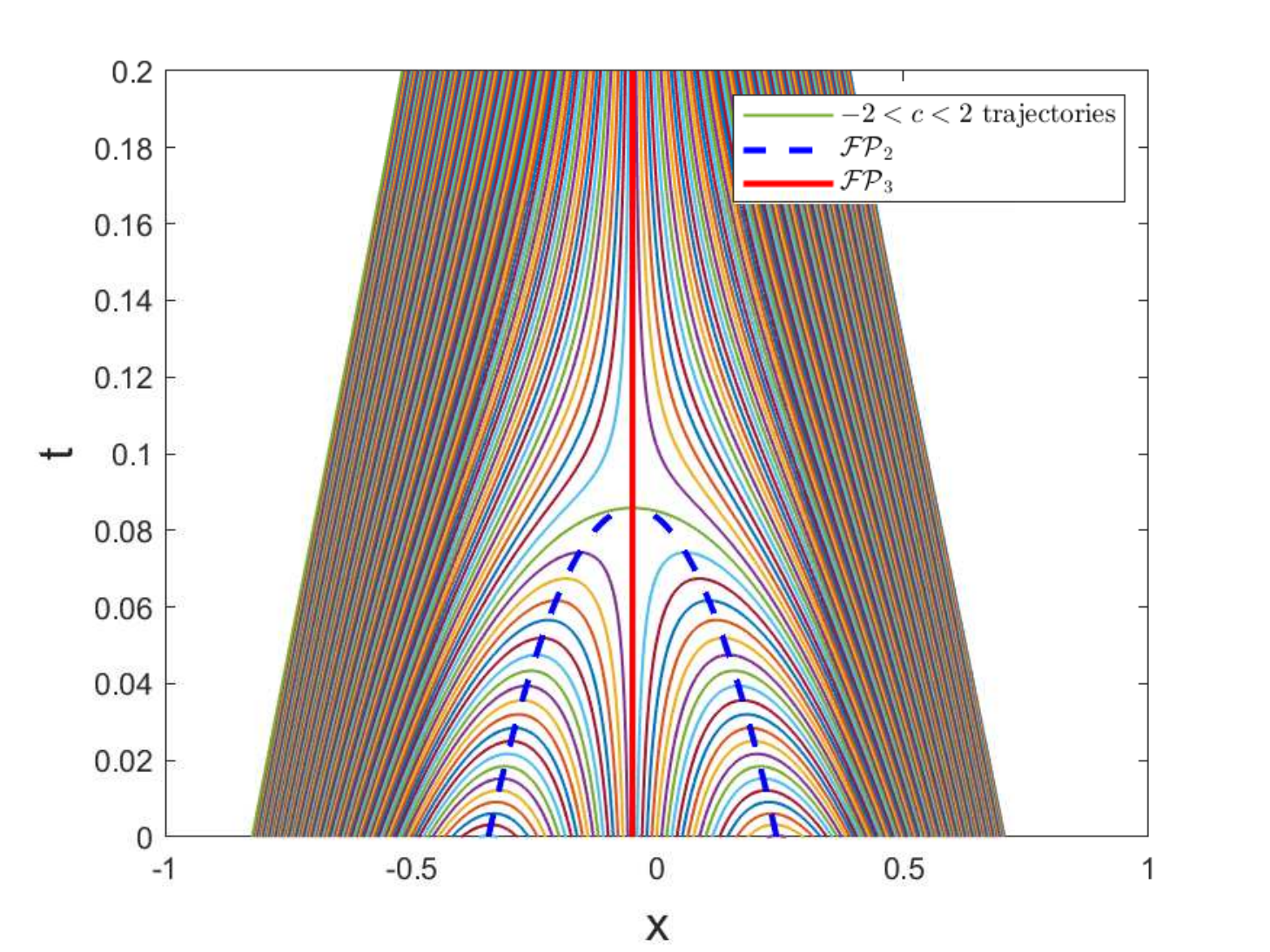}
\caption{}
\end{subfigure}
\caption{Illustration of Fingerprints and Trajectories in Example \ref{example-14}. To observe the change of trajectories with coefficient $c$ in the polynomial, we vary it in $c\in [-2,2]$. (a): $\mathcal{FP}_i,\ (i=1,2,3)$,  $t_u$ and $t^*$. Notice that $\mathcal{FP}_1$ corresponds to $c=0.01$. (b): $\mathcal{FP}_2,\ \mathcal{FP}_3$ and trajectories of $c=-0.05$. (c): critical  trajectories when $c=-0.051$, which are symmetric about $x=-\frac{a}{4}$.  (d): $\mathcal{FP}_2,\ \mathcal{FP}_3$ and trajectories of $c=-0.2$. (e): trajectories by varying $c$. (f): more trajectories by varying $c$. }
\label{fig-example-14}
\end{figure}


\begin{remark}
If there exists a critical point $x'$ at  $t=0$ such that $x'=-\frac{a}{4}$, then  \eqref{eqn:dyn-1d} implies its fingerprint curve $x(t)\equiv -\frac{a}{4}$. This happens when $x'$ is the local maximizer, and other two critical points $x_1$ and $x_3$ satisfy $x_1+x_3=2x'$. If case is this, all three
fingerprint curves meet up at $x'=-\frac{a}{4}$ when $t=t_u$,  which  will be explained in details in the following sections.
\end{remark}

\begin{example}\label{example-14} Consider the polynomial $p(x)=x^4+0.2x^3-0.5x^2+0.01x$,  the illustration is in Fig.\ref{fig-example-14}.
\end{example}

\subsection{Heat evolution of critical points of quartic polynomials}\label{Good}

Now we investigate the evolution of the critical points, and in particular, the quart polynomial case. Essentially, we concentrate on the case $t_u>0$ in which there exist three distinct critical points $ x_1<x_2< x_3$, and they correspondingly evolve to the critical points $x_1^t,x_2^t,x_3^t$ when $0<t<t_u$.
Generally, both $x_1$ and $x_3$ are local minimizers and $x_2$ is local maximizer. Our main concern is  the behavior associate with heat evolution,  characterized by equation \eqref{eqn:dyn-1d}.


Our next concern is the comparison principle between two local minimums during heat evolution. Surprisingly, we have the very expected  result for heat evolution of quartic polynomials: 

\begin{theorem}\label{main} For monic quartic polynomial $p(x)$, assume that $t_u>0$, and denote its three critical points $ x_1< x_2<x_3$ (or $x_1>x_2>x_3$). 
If $p(x_1)<p(x_3)$, then $p(x_1^t,t)<p(x_2^t,t)$.
\end{theorem}

Before showing this result, we need several lemmas.

\begin{lemma}\label{temple-1} Let $x_1,x_2,x_3$ be the critical points of quartic polynomial $p(x)$, then we have
\begin{equation}
p(x_3)-p(x_1)=-(x_3-x_1)^3\cdot (x_1+x_3-2x_2)/3.
\end{equation}
\end{lemma}
\begin{proof}
Represent $a,b,c$ in terms of $x_1,x_2,x_3$, by \eqref{abc-by-x123}. Thus we obtain $p(x_1)-p(x_3)$. Factorizing it will lead to required result.
\end{proof}

This Lemma \ref{temple-1} implies that
\begin{lemma}\label{lem-equal-peaks}
For any $t\in [0,t_u)$, $ p(x_3^t,t)=p(x_1^t,t)$ if and only if $ x_1^t+x_3^t=2x_2^t$.
\end{lemma}

Consequently, we see that
\begin{lemma}\label{lem-equal-peaks-1}
For any $t\in [0,t_u)$, $  x_1^t+x_3^t=2x_2^t$ if and only if $ x_2^t=-\frac{a}{4}$.
\end{lemma}
\begin{proof}
Notice that the coefficient of $x^3$ in $p(x,t)$ is invariant with $t$, and the coefficient of $x^2$ of $\frac{\partial p}{\partial x}$ is also invariant with $t$. According to Appendix \ref{appendix-1}, we have
\begin{equation} x_1^t+x_2^t+x^t_3=-\frac{3a}{4}.\end{equation}
Applying Lemma \ref{lem-equal-peaks} will yield the result.
\end{proof}

\begin{lemma}\label{equal-is-impossible} Assume that $x_1<x_2<x_3$  are three critical points of monic quartic polynomial $p(x)$, if $p(x_1)=p(x_3)$, then for all $t\in (-\infty,t_u)$, we have
\begin{equation}
p(x_1^t,t)=p(x_3^t,t).
\end{equation}
\end{lemma}
\begin{proof}
Recall \eqref{eqn:fp-1c}, and apply \eqref{abc-by-x123},  we can actually represent
$\frac{\partial p(x,t)}{\partial x}=0$  in terms of $x_1,x_2,x_3$ as
\begin{equation}
\begin{split}
x^3-(x_1+x_2+x_3)x^2+&(x_1x_2+x_2x_3+x_3x_1+3t)x\\
                 -&[x_1x_2x_3+(x_1+x_2+x_3)t]=0.
\end{split}
\end{equation}
Now if $p(x_1)=p(x_3)$, Lemma \ref{lem-equal-peaks} tells us $x_3=2x_2-x_1$, thus we may simplify the above equation as
\begin{equation}
x^3-3x_2x^2+(2x_1x_2+2x_2^2-x_1^2+3t)x-(2x_1x_2^2 -x_1^2x_2+3x_2t)=0,
\end{equation}
whose solution is
\begin{equation}
\begin{split}
x_2^t&=x_2,\\
x_1^t&=x_2-\sqrt{(x_1-x_2)^2-3t},\\
x_3^t&=x_2+\sqrt{(x_1-x_2)^2-3t},\ \ \ \ \ -\infty<t<\min\left \{\frac{(x_1-x_2)^2}{3},t_u\right \}
\end{split}
\end{equation}
which shows that $x_1^t+x_3^t=2x_2^t$. According to Lemma \ref{lem-equal-peaks}, this leads to  $p(x_1^t,t)=p(x_3^t,t)$.
\end{proof}

At present stage, we summarize all lemmas as below,

\begin{theorem}\label{main-again} Under the same assumptions as Theorem \ref{main}, and denote $x_1^t< x_2^t< x_3^t$ the critical points of $p(x,t)$. Then the following statements are equivalent:
\begin{enumerate}
\item $
p( x_1)=p( x_3),
$
\item $
p( x_1^t,t)=p( x_3^t,t),\ \forall t\in [0,t_u )
$;
\item $x_1+x_3=2x_2$;
\item $x_1^t+x_3^t=2x_2^t,\ \forall t\in [0,t_u)$;
\item $x_2=-\frac{a}{4}$;
\item $x_2^t=-\frac{a}{4},\ \forall t\in [0,t_u)$.
\item $t^*=t_u$.
\end{enumerate}
\end{theorem}

\begin{proof}
We will prove that (1)$\implies$(3)$\implies$(5)$\implies$(6)$\implies$(4)$\implies$(2)$\implies$(1). In addition, (4)$\iff$(7). Actually, this routine is partially repeated with previous proofs.

(1)$\implies$(3) (also (4)$\implies$(2)) comes from Lemma \ref{lem-equal-peaks},  (3)$\implies$(5) (also (6)$\implies$(4)) from Lemma \ref{lem-equal-peaks-1}, (5)$\implies$(6) from differential equation \eqref{eqn:dyn-1d}. Finally, (4)$\iff$(7) comes from \eqref{t-start-t-u} in Lemma \ref{13} as well as the condition $x_1<x_2<x_3$.
\end{proof}

\begin{proof}[Theorem \ref{lem-equal-peaks}]
The dynamical equation \eqref{eqn:dyn-1d} states that the  evolution of three critical points are continuous when $t\in [0,t_u)$. Thus if $p(x_1)<p(x_3)$,  we must have $p(x_1^t,t)<p(x_3^t,t)$ for $t\in [0,t_u)$, otherwise, there must have
$$p(x_1^{t'},t')=p(x_3^{t'},t')$$
for some $t' \in (0, t_u)$. However, if this is the case, Theorem \ref{main-again} (Lemma \ref{equal-is-impossible}) indicates that $p(x_1) = p(x_3)$ because $t$ is reversible, which contradicts the initial assumption.
\end{proof}

To intuitively explain this result, we suggest a triangle representation at Fig.\ref{Triangle of quartic polynomial} for each $t$, where the cortes of triangle consists of $(x_i^t,p(x_i^t,t)),\ (i=1,2,3)$ when $t<t_u$. Notice that the sequence of triangles when $0\le t\le t_u$ and continued curve $\dot x(t)$ actually connected to global minimum of $p(x,t)$ at each $t\ge 0$.

Finally, we discuss an interesting problem: if $x_1<x_2<x_3$ are three critical points of quartic polynomial $p(x)$, can we judge which one of them is global minimizer without valuating all these $p(x_i)$? The answer is YES.

\begin{theorem}\label{th:A} Let $x_1<x_2<x_3$ be three distinct critical points of monic quartic polynomial $p(x)$, then the following statements are equivalent:
\begin{enumerate}
\item $x_3$ (resp. $x_1$) is global minimizer;
\item $x_1+x_3>2x_2$\ (resp. $x_1+x_3<2x_2$);
\item  $x_2<-a/4$ (resp. $x_2>-a/4$).
\end{enumerate}
\end{theorem}
\begin{proof} Apply Lemma \ref{temple-1}.
\end{proof}
This Theorem \ref{th:A} inspired the following very simple Euler's algorithm without Heat Convolution for quartic polynomials.
\begin{theorem} For any monic quartic polynomial $p(x)$ with $a$ the  coefficient of $x^3$, the Euler's algorithm  with FIXED initial position $x^{(0)}=-\frac{a}{4}$,
\begin{equation}
x^{(k+1)}=x^{(k)}-\Delta x\cdot p'(x^{(k)}),
\end{equation}
MUST converge to the global minimizer of $p(x)$.

\end{theorem}

\subsection{Algorithm}
 Actually, the Euler's algorithm may work from $t>t_u$. Fortunately, we know at $t\ge t_u$, the $p(x,t)$ has only single minimum about $x$. Recall
the formula \eqref{appendix-2}, in which we see that the sum of all real roots of fingerprint cubic equation \eqref{eqn:fp-1c} should be invariant under $0\le t\ge t_u$, thus we know the remaining critical point $x_{init}$ at $t=t_u$ can be solved since we already know the information of $(x(t_u),t_u)$ from Theorem \ref{singularity-xt}. This means we may adopt
\begin{equation}
\begin{split}
x_{init}=&-\frac{3a}{4}-2x(t_u)\\
=&-\frac{a}{4}-2\left (\frac{a^3-4ab+8c}{64} \right )^{1/3},\\
\end{split}
\end{equation}
at $t=t_u$ as initial point, then perform Euler's algorithm for equation \eqref{eqn:dyn-1d}, and finally attain the global minimum of $p(x)$.
This implies the following result.

\begin{theorem} For the quartic polynomial \eqref{a1}, if $t_u\le 0$, this polynomial has only one critical point. If $t_u>0$, the polynomial contains three distinct critical points, then if all the critical points satisfy $x\ne -\frac{a}{4}$, we backward perform the differential equation
\begin{equation}\label{eqn:dyn-1a}
\frac{\,d{ x}(t)}{\,dt}=-\frac{12{x}+3a}{12{x}^2+6ax+2b+12t},
\end{equation}
with initial condition
\begin{equation}
x^{t_u}=x(t_u)=-\frac{a}{4}-\frac{(a^3-4ab+8c)^{1/3}}{2}
\end{equation}
 from $t=t_u>0$ to $t=0$, must attain  the global minimizer of \eqref{a1} at $t=0$. Finally, if one critical point equals $-\frac{a}{4}$, then $p(x)$ has two
 global minimizers, and they are the roots of quadratic polynomial
 \begin{equation}
 \frac{p(x)}{x+\frac{a}{4}}.
 \end{equation}
\end{theorem}

So far, we may start with verifying that whether $-\frac{a}{4}$ is a root of cubic polynomial $\frac{\partial p(x)}{\partial x}$.  The global minimizer can be obtained immediately if $-\frac{a}{4}$ is a root. Otherwise, set $x^{(0)}=x_{init}$, and $t^{(0)}=t_u$. Then motivated by \eqref{eqn:dyn-1d}, the iteration process is as below
\begin{equation}
\begin{split}
   x^{(i+1)}=&x^{(i)}-\Delta t\cdot \frac{12{x^{(i)}}+3a}{12{x^{(i)}}^2+6ax^{(i)}+2b+12t^{(i)}},\\
   t^{(i+1)}=&t^{(i)}-\Delta t.\\
 \end{split}
\end{equation}
Here the prescribed step $\Delta t>0$  is small enough, and we may stop the iteration while  $t^{(n)}\approx 0$. Finally, this algorithm provides
\begin{equation}
\lim_{i\to\infty} x^{(i)}=x_{min}.
\end{equation}

Instead beginning with $x_{init}$, we can also start by sufficient evolution $p(x,t)$. This will cost more steps of iterations.

\subsection{Numerical experiments}

\begin{landscape}
\begin{figure}
  \centering
  \includegraphics[width=22cm]{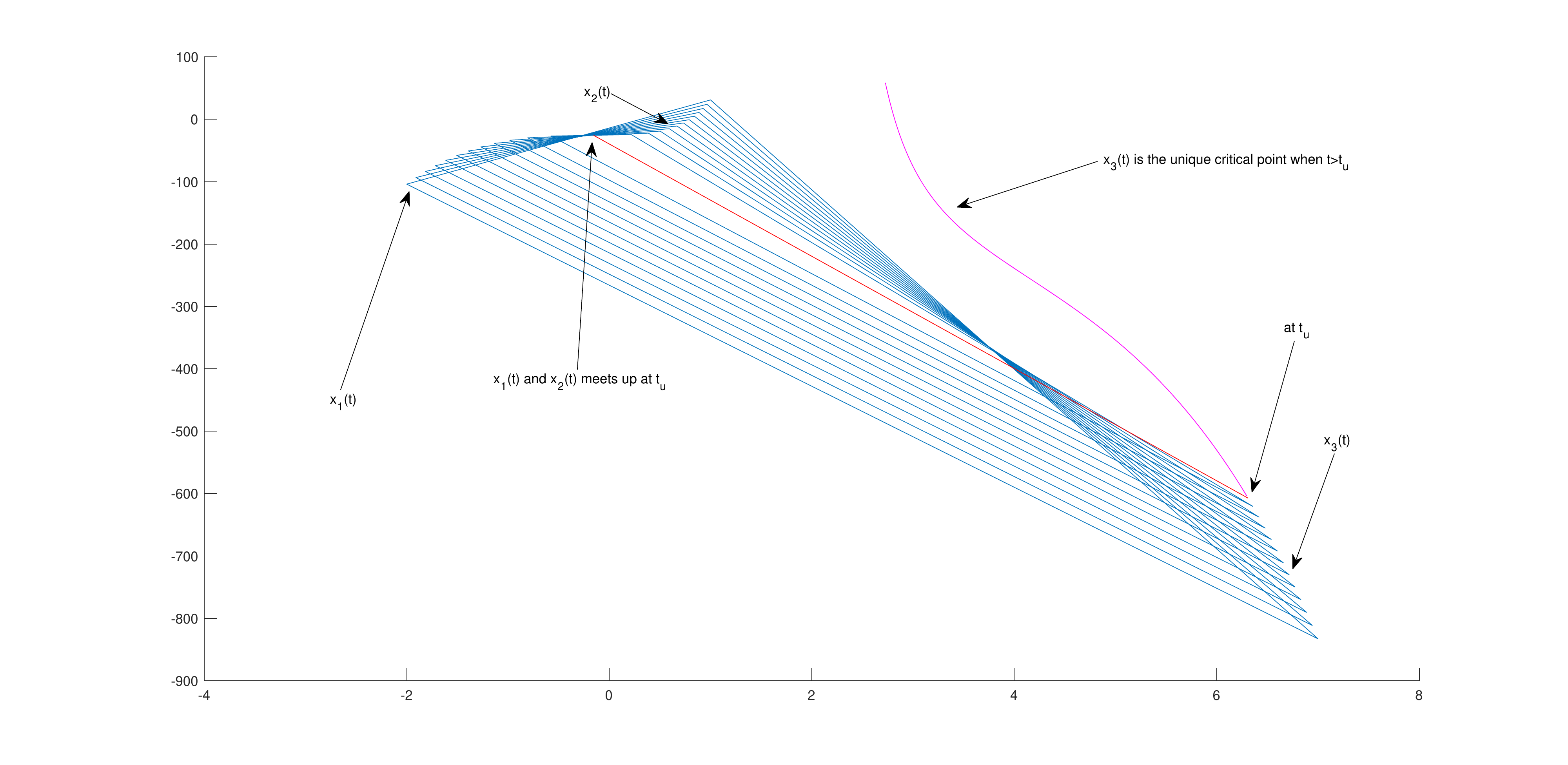}\\
  \caption{The triangle series of evolution of critical points. Here  $p(x)=x^4-8x^3-18x^2+56x$ is analyzed in Example \ref{example-1}. }\label{Triangle of quartic polynomial}
\end{figure}
\end{landscape}

\begin{example}\label{example-1} This  counter-example \cite{ABK-1} is proposed against 'backward differential flow' method of \cite{Zhu-2}, in which $p(x)=x^4-8x^3-18x^2+56x$. In our heat conduct framework, we have $p(x,t)=x^4-8x^3-(18-6t)x^2+32x-(18t-3t^2)$.   Notice that Fig.\ref{Triangle of quartic polynomial} 
demonstrates the triangle series of critical points.
\end{example}

\begin{figure}[ht]
  \centering
  \includegraphics[width=14cm]{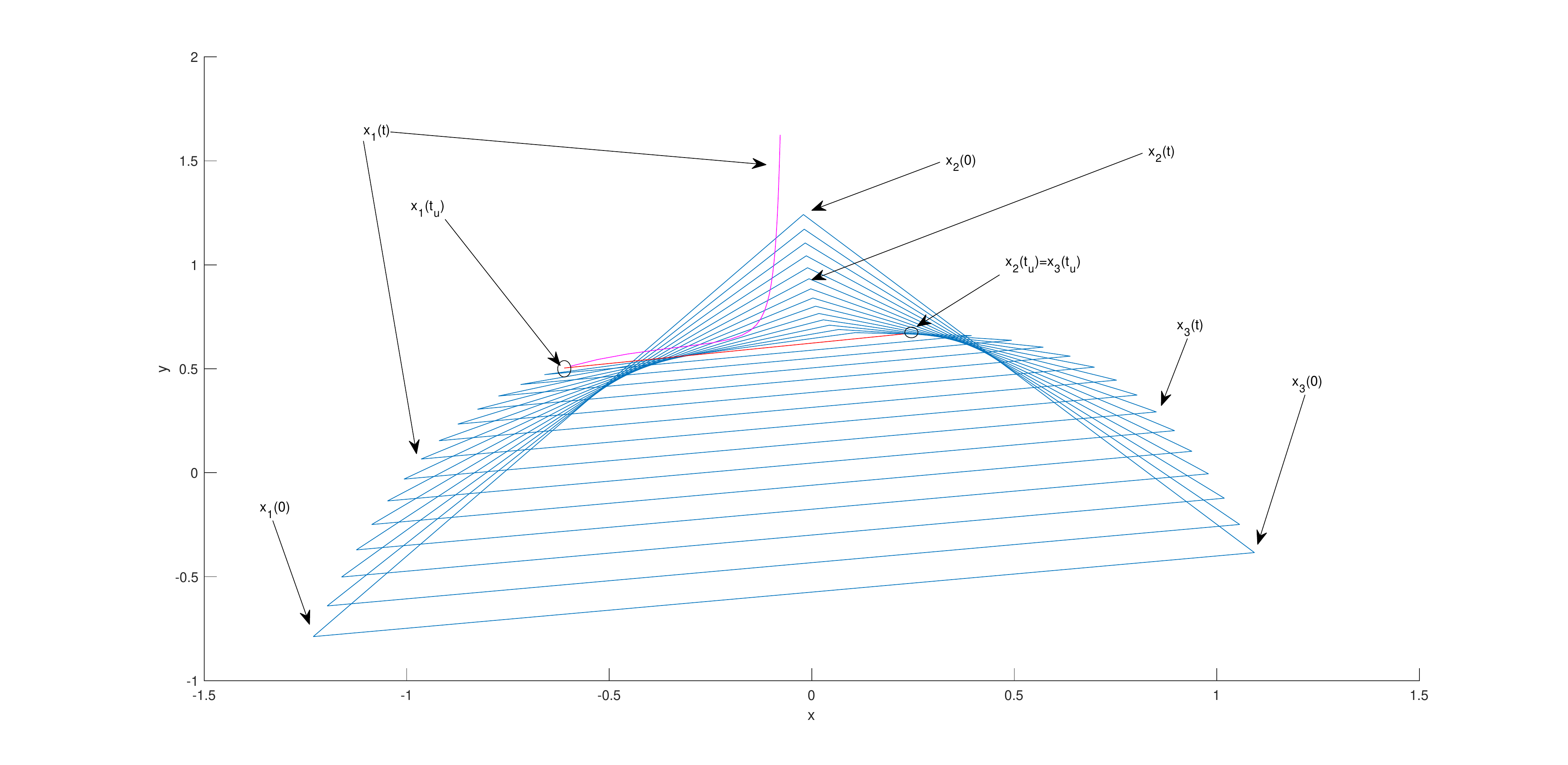}\\
  \caption{An example of triangle series in $(x,y)$  system, of $p(x)=x^4+0.2114 x^3  -2.6841x^2 -0.1110x+  1.2406$. See Example \ref{example-2} for more details.}\label{quart-02114}
\end{figure}

\begin{figure}[ht]
  \centering
  \includegraphics[width=8cm]{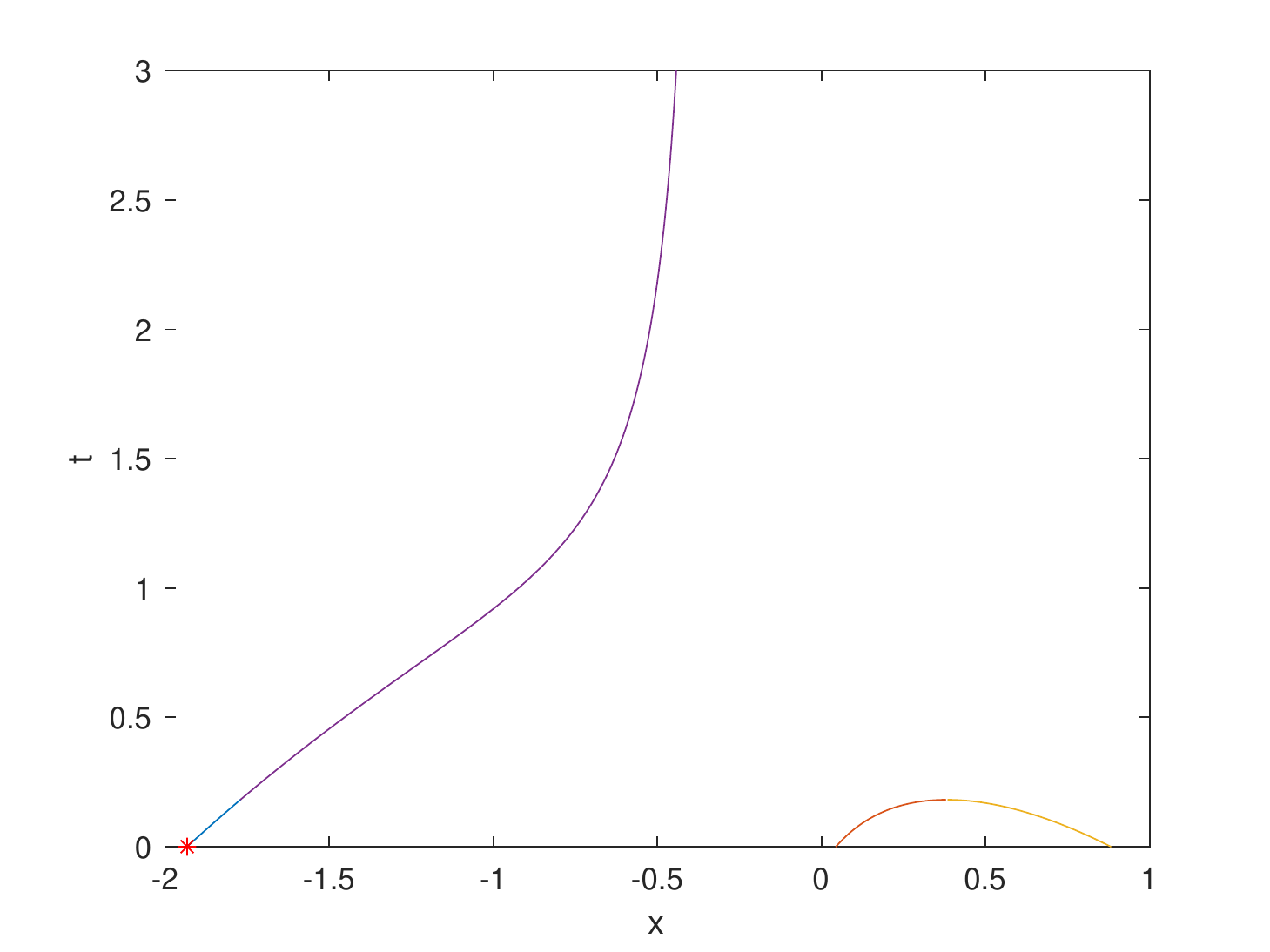}\\
  \caption{Fingerprint $\mathcal{FP}_1$ in $(x,t)$ system, of $p(x)=x^4+0.2114 x^3  -2.6841x^2 -0.1110x+  1.2406$. See Example \ref{example-2} for more details.}\label{quart-02114-fp}
\end{figure}

\begin{example}\label{example-2}Set $p(x)=x^4+0.2114 x^3  -2.6841x^2 -0.1110x+  1.2406$, then  in Fig. \ref{quart-02114}  the most left curve  $x_1(t)$ is the global minimizer of
corresponding $p(x,t)$ at each $t\ge 0$, and  Fig.\ref{quart-02114-fp} illustrates the fingerprint $\mathcal{FP}_1$. The theoretical minimizer is $x_1=-1.2307$,
and our iteration algorithm provides $x_1=-1.2308$.
\end{example}

\begin{example}\label{symmetric-ex} Consider $p(x)=x^4-4x^3-2x^2+12x$, then we actually have three critical points $x_1=-1,\ x_2=1,\ x_3=3$. Notice that
$x_1+x_3=2x_2$ thus $p(x_1,t)=p(x_3,t)$ and $x_2(t)=x_2=1$ for all $x\in [0,t_u]$. One can further verify that in this symmetric case, we must have $t^*=t_u$. The detailed explain can be referred as in Theorem \ref{main-again}.
\end{example}

\subsection{Summary of quartic polynomial case}
For the global minimizer of  quartic polynomial $p(x)$, while generating  its multi-scale version $p(x,t)=p(x)*g_t(x)$ on account of  Gaussian filter $g_t(x)$ with variance from $t=0$ to $+\infty$, we will see that:
\begin{itemize}
 \item If $t_u<0$, $p(x)$ itself is not necessary convex, but it has unique critical point. Consequently, each $p(x,t)$ has only one critical point at any $t\ge 0$;
 \item If further $t_u\le t^*<0$,  $p(x)$ must be convex. Then each $p(x,t)$ is convex about $x$ at any $t\ge 0$;
 \item If $t_u>0$, the polynomial $p(x)$ has three distinct critical points $x_1<x_2<x_3$ when $0\le t<t_u$;
\item When $t=t_u$, the critical point corresponding to the global minimizer will evolve continuously from $t_u$ to $t^*$, and the local minimizer will meet up with local maximum $x_2^t$. Even more, these two critical points will stop evolution at $t=t_u$;
\item When $t_u<t<t^*$, the polynomial $p(x,t)$ has unique minimizer at each $t$.
\item When $t\ge t^*$, the polynomial $p(x,t)$ will become  convex about $x$, and possesses unique minimizer.
\end{itemize}

\section{Case study  of sixth degree polynomials}
\subsection{Evolution and fingerprints}
Now we consider six degree monic polynomial
\begin{equation}\label{6-1}
p(x)=x^6+bx^4+cx^3+dx^2+ex+f,
\end{equation}
For the sake of simplicity,  here we regularize the coefficient of $x^5$ by setting it as zero, which is a standard technique in treating the algebraic equations. Thus the heat evolution is
\begin{equation}\label{6-2}
\begin{split}
p(x,t)=&p(x)+t\cdot\frac{\partial p}{\partial t}+\frac{t^2}{2}\frac{\partial^2 p}{\partial t^2}+\frac{t^3}{6}\frac{\partial^3 p}{\partial t^3}\\
     =&p(x)+\frac{t}{2}\cdot\frac{\partial^2 p}{\partial x^2}+\frac{t^2}{8}\frac{\partial^4 p}{\partial x^4}+\frac{t^3}{48}\frac{\partial^6 p}{\partial x^6}\\
     =& x^6+b(t)x^4+c(t)x^3+d(t)x^2+e(t)x+f(t),\\
     \end{split}
\end{equation}
in which
\begin{equation}\label{6-3}
\begin{split}
b(t)&=b+15t,\\
c(t)&=c,\\
d(t)&=d+6bt+45t^2,\\
e(t)&=e+3ct,\\
f(t)&=f+dt+3bt^2+15t^3.\\
\end{split}
\end{equation}

The critical points of \eqref{6-1} satisfy the 5 degree equation
\begin{equation}\label{6-4}
0=\frac{1}{6}\frac{\partial p(x,t)}{\partial x}=x^5+B(t)x^3+C(t)x^2+D(t)x+E(t),
\end{equation}
in which
\begin{equation}\label{6-5}
\begin{split}
B(t)&=\frac{2b}{3}+10t,\\
C(t)&=\frac{c}{2},\\
D(t)&=\frac{d}{3}+2bt+15t^2,\\
E(t)&=\frac{e}{6}+\frac{ct}{2}.
\end{split}
\end{equation}
The solution of \eqref{6-4} for $t\ge 0$ consists the fingerprint $\mathcal{FP}_1$.
Unfortunately, the roots of this fifth degree equation is algebraically intractable \cite{Stewart}.

Similarly, we may write the equation of $\mathcal{FP}_2$ as below,
\begin{equation}\label{6-6aaa}
\frac{1}{30}\frac{\partial^2 p }{\partial x^2}=x^4+\left (\frac{2b}{5}+6t \right  ) x^2+\frac{c}{5}x+\frac{d}{15}+\frac{2b}{5}t+3t^2 =0.
\end{equation}

Then  the equation of $\mathcal{FP}_3$ is
\begin{equation}\label{6-6ab}
\frac{1}{120}\frac{\partial^3 p }{\partial x^3}=x^3+ \left (\frac{b}{5}+3t\right )x+\frac{c}{20}=0.
\end{equation}

Our interest is the intersection $\mathcal{FP}_2\bigcap\mathcal{FP}_3$, which offers the merge information.  Different from quartic polynomials case, we have not any explicit representation right now. Clearly, the real double root of quartic equation \eqref{6-6aaa} must be the common roots of both \eqref{6-6aaa} and \eqref{6-6ab}. In general, there exist two $0\le t_1<t_2$ such that the corresponding $x_1$ and $x_2$ are those two real double roots. The following Theorem \ref{Thm-t-x} explains the process of numerical approach.



\begin{figure}[t]
  \centering
  \includegraphics[width=10cm]{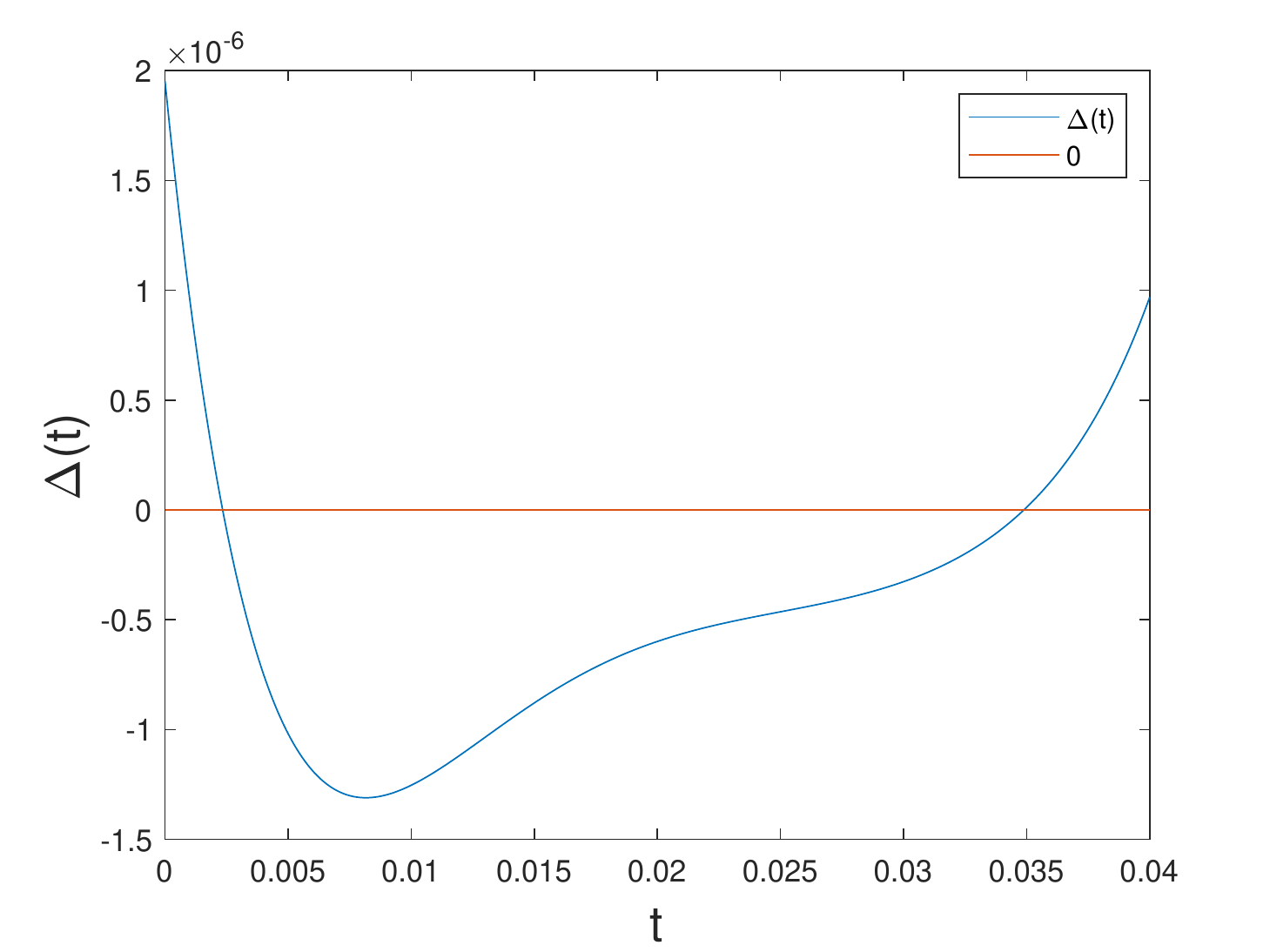}\\
  \caption{The discriminant $\Delta(t)$ of quartic equation is generated from $\frac{\partial^2 p(x,t) }{\partial x^2}=0$ and defined in \eqref{Thm-t-x-1}. Here the data comes from Example \ref{example-0}. The two real roots of $\Delta(t)$ are $t_1=0.002341$, at which the quartic equation possesses a real double root and two distinct real roots, and $t_2=0.034887$, at which the quartic equation possesses a real double root and a pair of conjugate complex roots.}\label{discrimant_six}
\end{figure}

\begin{figure}[b]
  \centering
  \includegraphics[width=10cm]{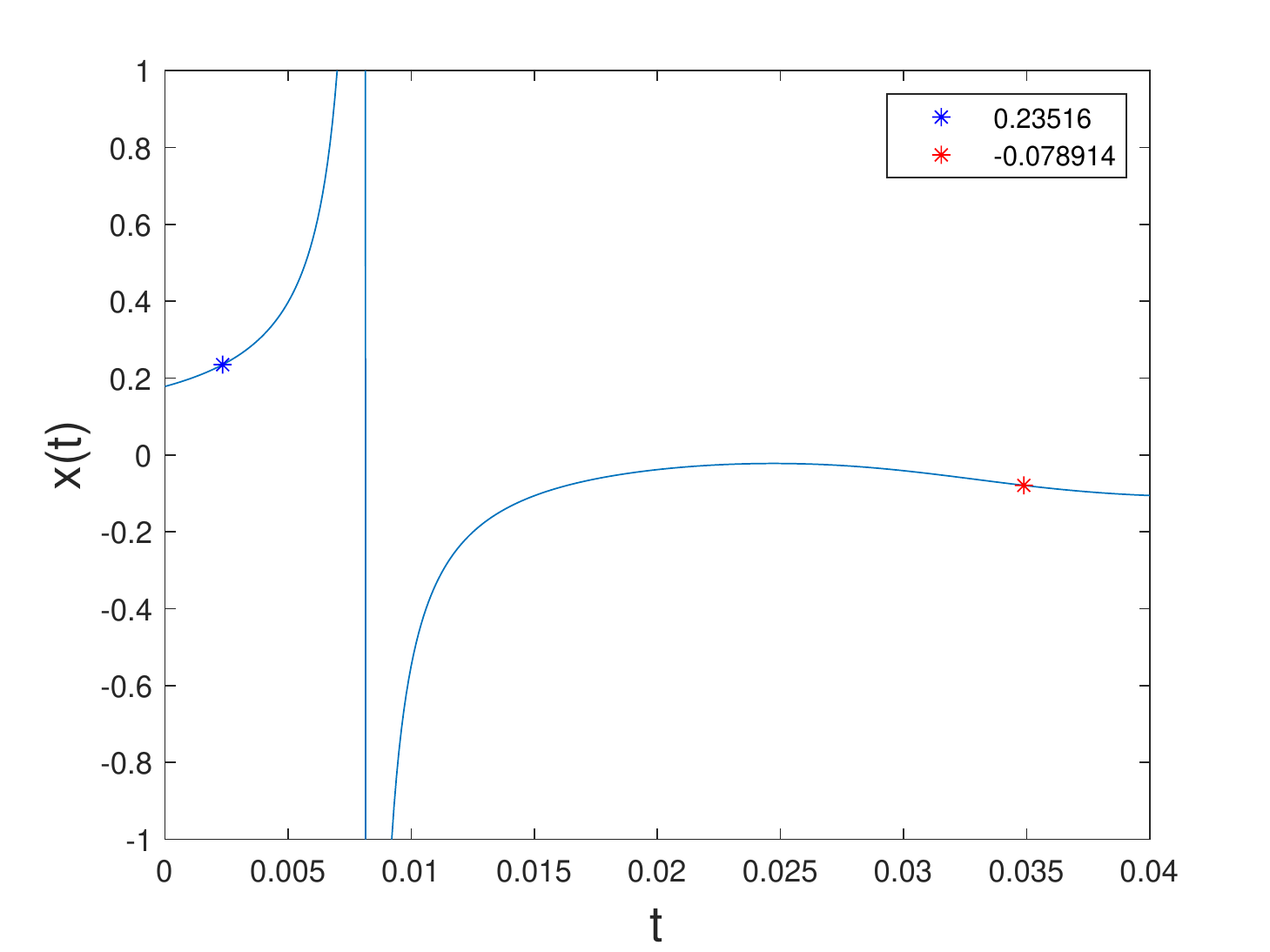}\\
  \caption{The function $x(t)$ is defined in \eqref{6-6ae}, where the data comes from the Example \ref{example-0}. Here we obtain two solution $(t_1,x_1)=(0.0023,   0.23516)$, $(t_2,x_2)=  (0.03489,-0.078914)$, which consist the set $\mathcal{FP}_2\bigcap\mathcal{FP}_3$. }\label{x_of_t}
\end{figure}

\begin{theorem}\label{Thm-t-x} For any six degree monic polynomial $p(x)$, the set $\mathcal{FP}_2\bigcap\mathcal{FP}_3$ contains a pair of elements $(x_i,t_i),\ i=1,2$, or
one element $(x_1,t_1)$, or  empty. Specifically, \newline
(1) Any $t_i$ must be the zero of discriminant
 \begin{equation}\label{Thm-t-x-1}
  \Delta(t)=27648 t'^6+ \frac{1728c^2}{25} t'^3 -\frac{256}{625} h^2   t'^2 - \frac{288}{625}c^2h  t' -\frac{256h^3}{75^3}-\frac{27c^4}{625}
   \end{equation}
 where
 \begin{equation}
  t'= t+\frac{b}{15},\ \ \ h=b^2-5d.
 \end{equation}
(2) Any $x_i$ is dependent of $t_i$ by the function
\begin{equation}\label{6-6ae}
x = -c\cdot \frac{1800\left (t+\frac{b}{15} \right )^2-4(b^2-5d)}{36000\left (t+\frac{b}{15} \right )^3+  80(b^2-5d) \left (t+\frac{b}{15} \right )+45c^2   }.
\end{equation}
\end{theorem}

\begin{proof}
We will give a detailed analysis in Apendix \ref{solution}, based on which, we know that in general settings there exist at most  two {\it{merge time}} $t_1$ and $t_2$ from the discriminant equation $\Delta(t)=0$ of quartic equation \eqref{6-6aaa} and \eqref{App-B-6-delta}.  This function \eqref{Thm-t-x-1} can be verified by \eqref{App-B-delta} and \eqref{App-B-6-delta} immediately, but we omit the detailed computation here. 

Now we may take Euclidean algorithm to reduce the degree of the polynomials about $x$ for \eqref{6-6aaa} and \eqref{6-6ab}. At first,  multiplying \eqref{6-6ab} with $x$, and subtracted both sides of \eqref{6-6aaa} resp., we get a second degree polynomial
\begin{equation}\label{6-6ac}
\left (t+\frac{b}{15} \right )x^2+\frac{c}{20}x+\left (t+\frac{b}{15}\right )^2-\frac{b^2-5d}{225}=0.
\end{equation}
Again, multiplying with $x$ for both sides of \eqref{6-6ac}, and subtracted from  both sides of \eqref{6-6ab} multiplied  with $t+\frac{b}{15}$, then we obtain
\begin{equation}\label{6-6ad}
-\frac{c}{20}x^2+ \left [ 2\left ( t+\frac{b}{15} \right )^2+\frac{b^2-5d}{225} \right ] x+\frac{c}{20}\left ( t+\frac{b}{15}\right )=0.
\end{equation}

Finally,  eliminating the second degree term by combining \eqref{6-6ad} and \eqref{6-6ac} will lead to \eqref{6-6ae}.
\end{proof}

As explained in Appendix \ref{solution}, we can obtain the suitable $t$ from the discriminant $\Delta(t)=0$ at first, then substitute it into the above
\eqref{6-6ae}, then choose $x$ from \eqref{6-6ae}.

\begin{figure}
  \centering
  \includegraphics[width=13cm]{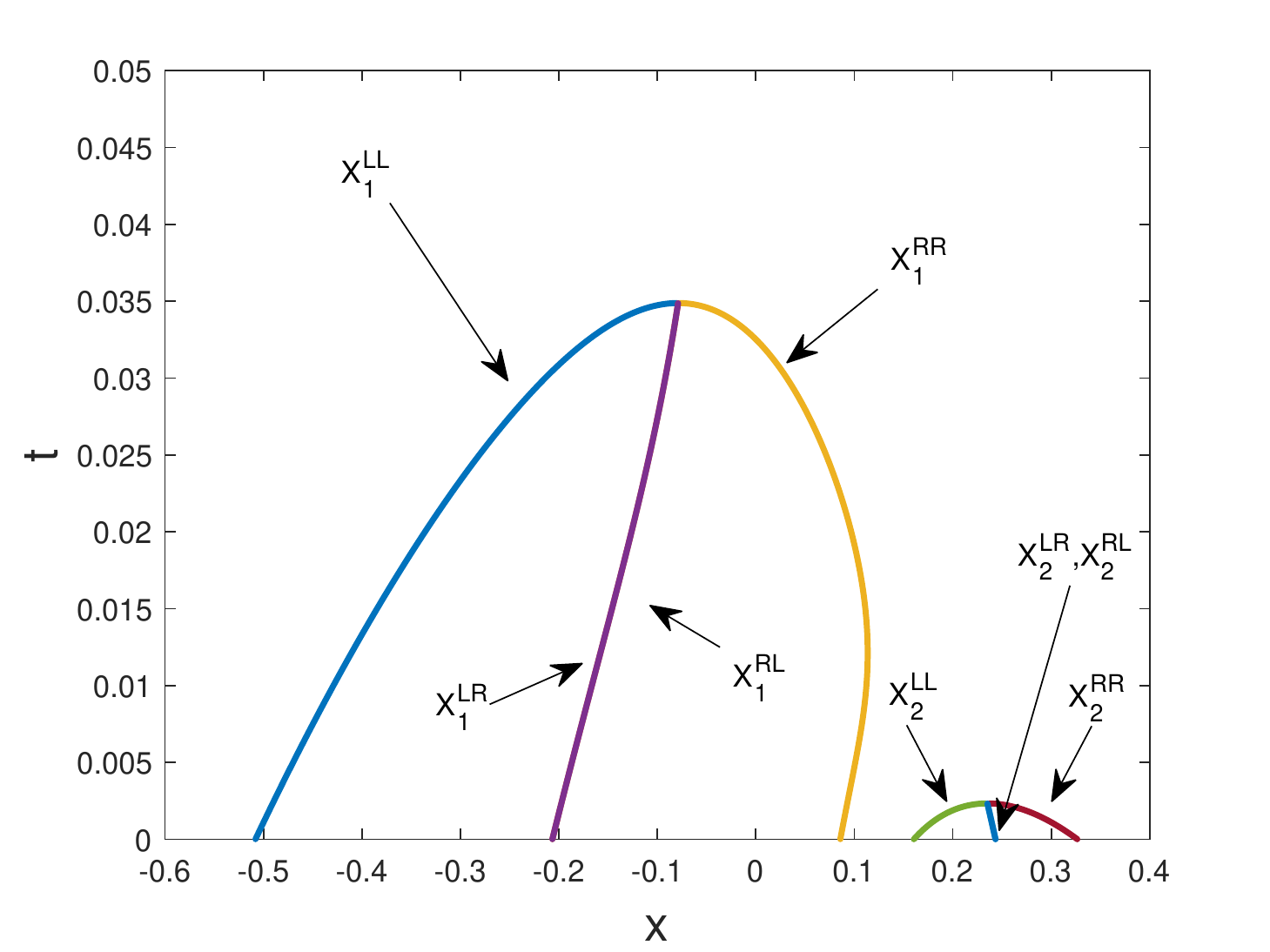}\\
  \caption{  Numerically, the partition of {\it Confinement Zone} and {\it Escape Zone} associated with $p(x)$ defined in Example \ref{example-0} is obtained through Matlab ODE packet of @ode25, and the {\it Confinement Zone}  is $[-0.5082,    0.0858]\bigcup [    0.1603,    0.3267]$.}\label{escape_zone}
\end{figure}

\subsection{Numerical examples for 6-degree polynomials}
It is extremely expected to generalize the heat evolution algorithm to find out global minimizer of 6 or higher even degree polynomials. However, the Theorem \ref{main} can not be generalized to higher degree polynomials. Here we illustrate the positive and negative examples.

\begin{figure}
\centering

\begin{subfigure}[a]{0.70\linewidth}
\includegraphics[width=\linewidth]{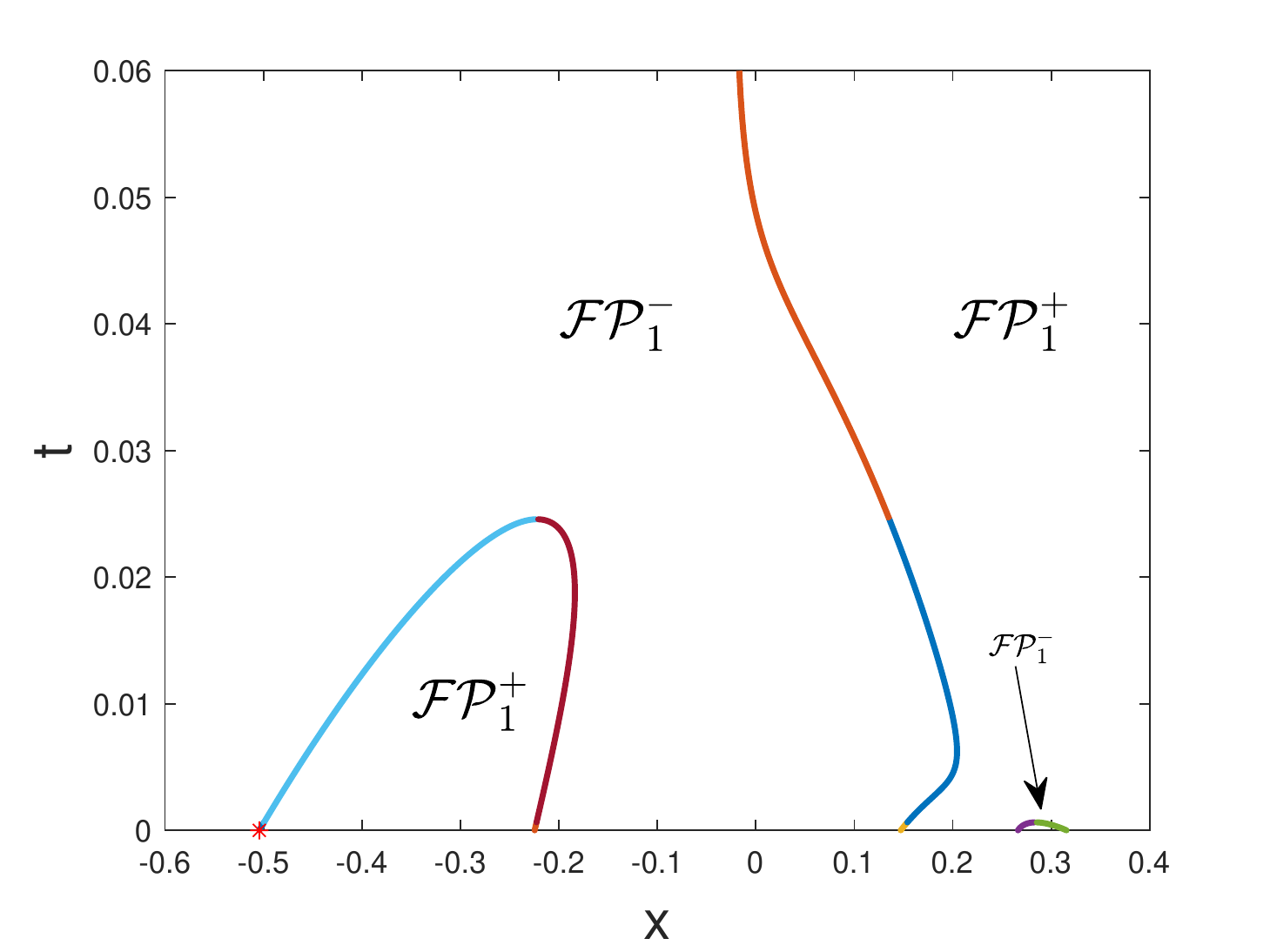}
\end{subfigure}

\begin{subfigure}[b]{0.70\linewidth}
\includegraphics[width=\linewidth]{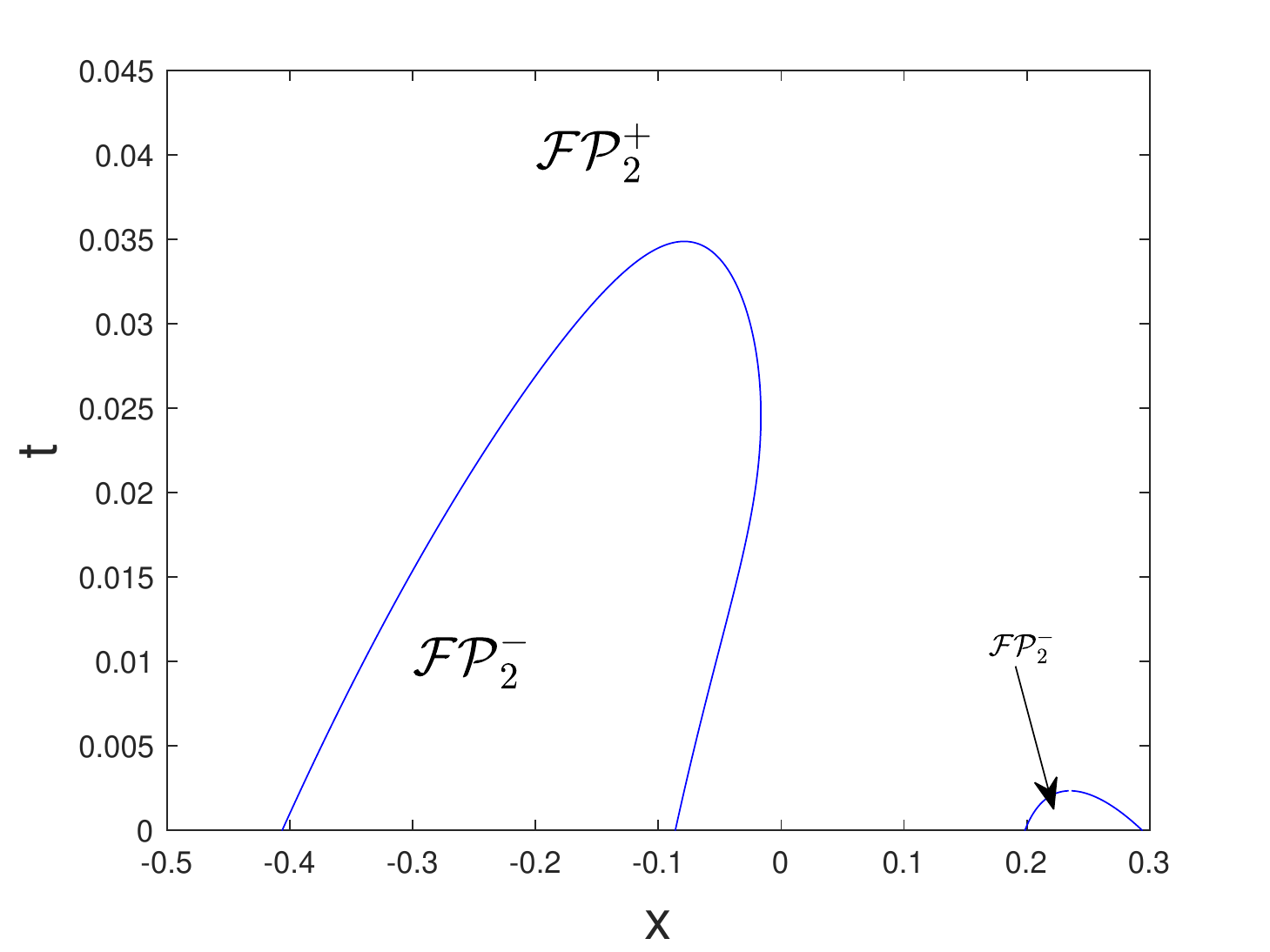}
\end{subfigure}

\begin{subfigure}[c]{0.70\linewidth}
\includegraphics[width=\linewidth]{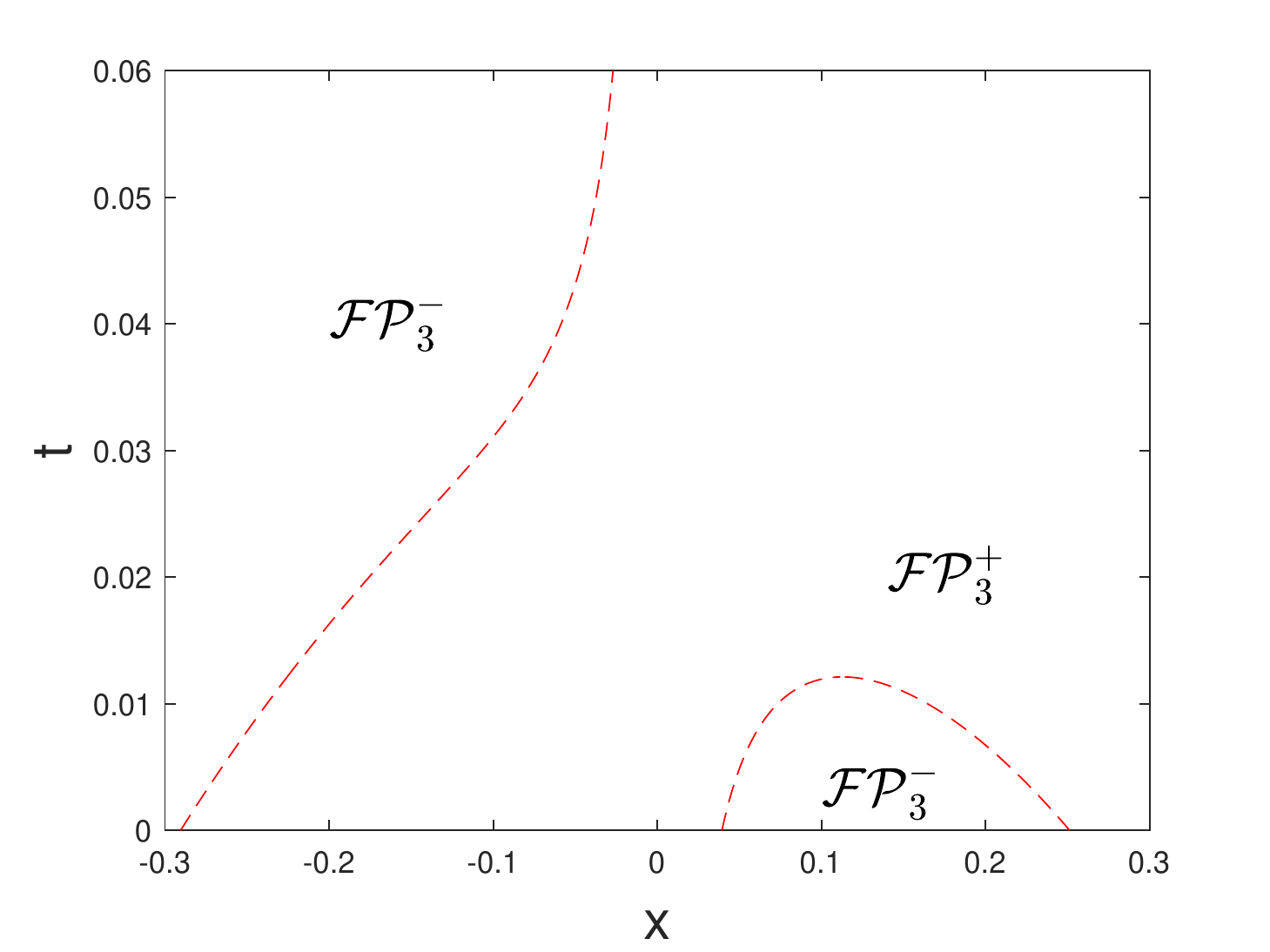}
\end{subfigure}

\caption{The fingerprints in Example \ref{example-6d-1}, separately illustrated.}
\label{fig-11A} 
\end{figure}

\begin{figure}{
\includegraphics[width=0.70\linewidth]{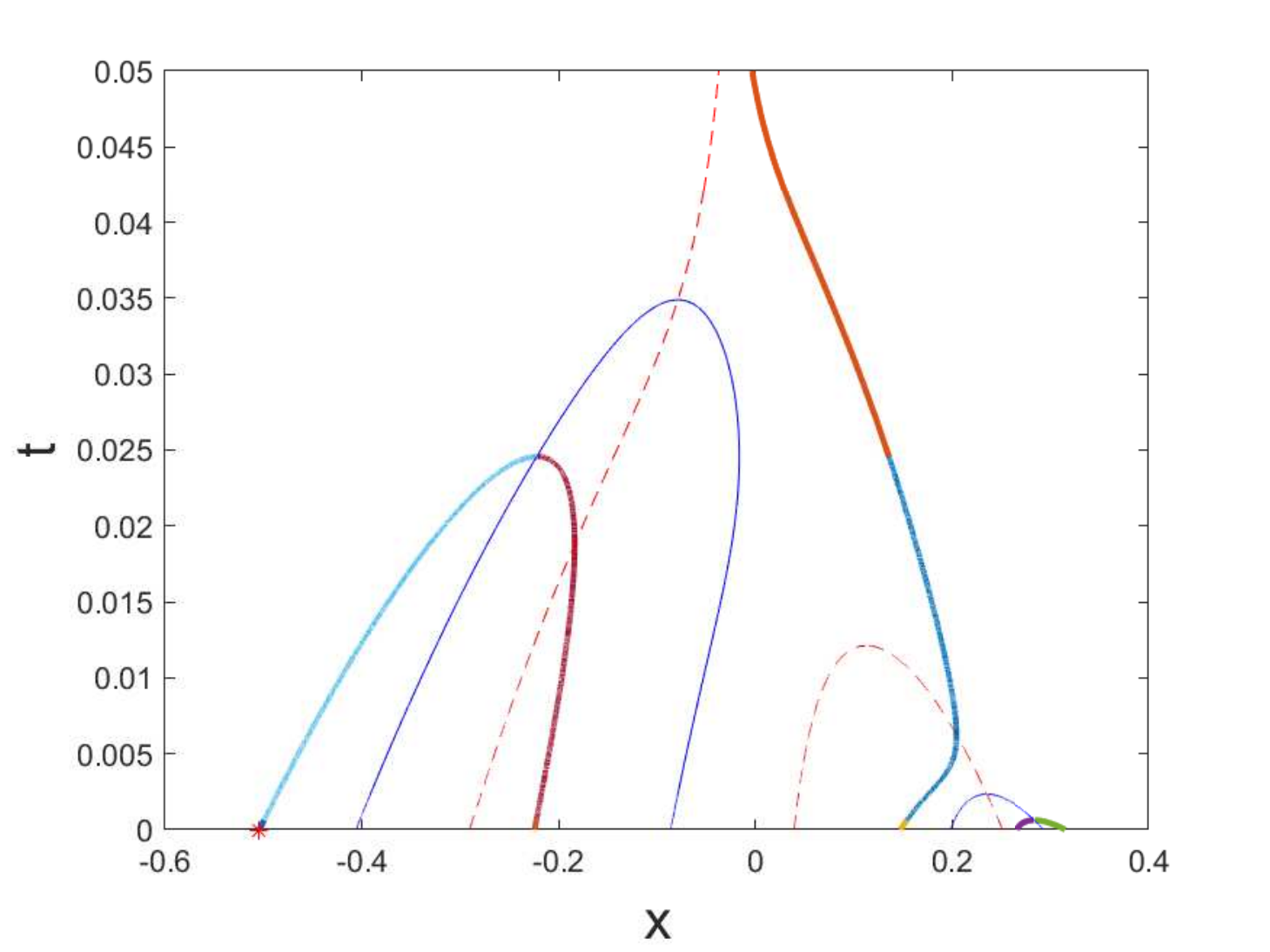}}
\caption{The joint illustration of fingerprints $\mathcal{FP}_1$, $\mathcal{FP}_2$ and $\mathcal{FP}_3$ of previous figures about Example \ref{example-6d-1}.}
\label{fig-m}
\end{figure}

\begin{example}\label{example-6d-1} We illustrate the fingerprints of six degree polynomial
$$ p(x)=x^6-0.3726x^4+0.0574x^3+0.0306x^2-0.0084x$$
in Fig.\ref{fig-11A} and Fig.\ref{fig-m}. We point out  that the global minimizer ("*" in Fig.\ref{fig-11A}(a)) does not evolute to large $t$, which means the convex convolution for this $p(x)$ will not converge to its global minimizer.
\label{example-0}
\end{example}

\begin{figure}
  \centering
  \includegraphics[width=9cm]{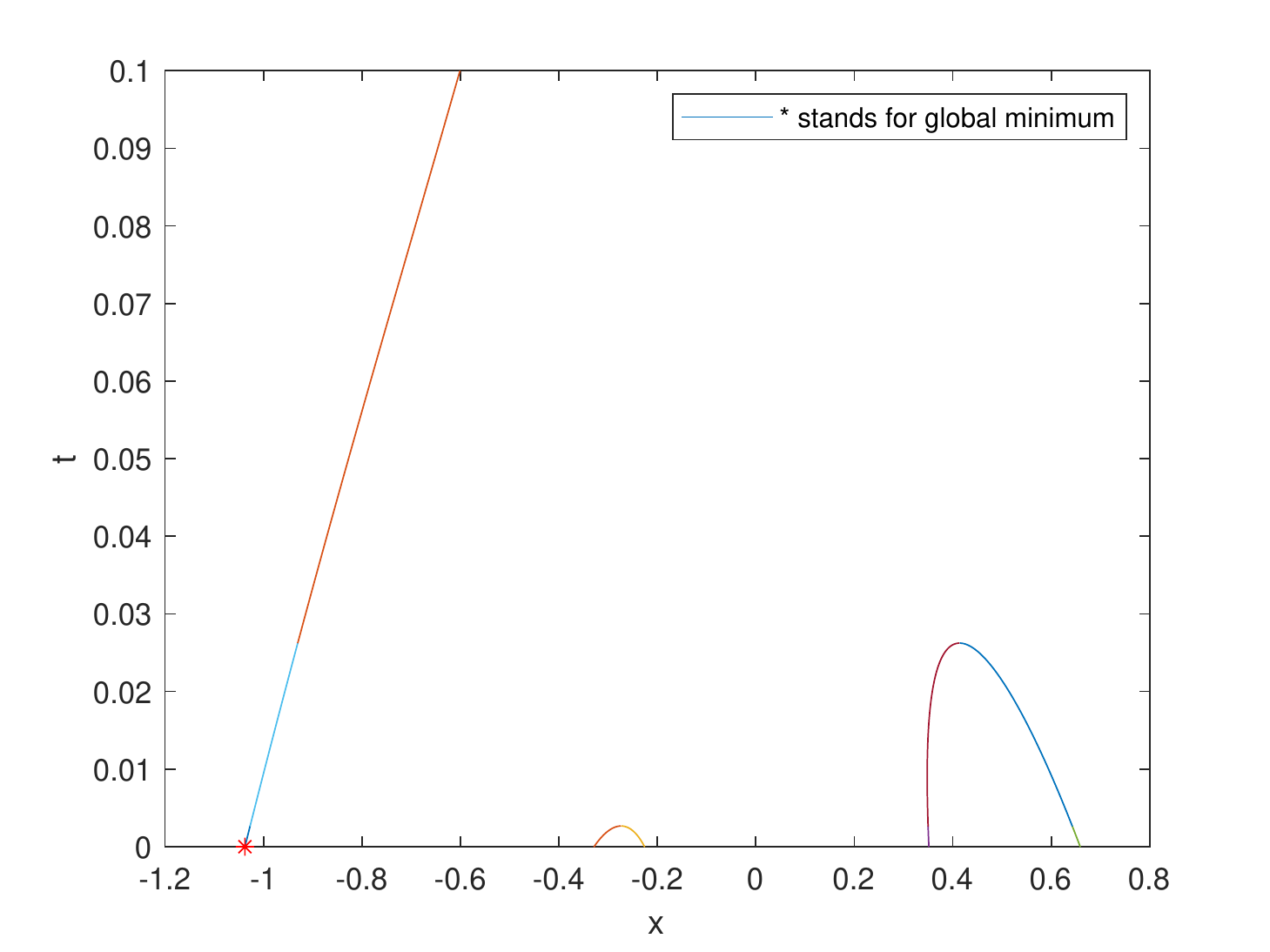}\\
  \caption{An example of Fingerprint of $p(x)=x^6+0.6987x^5-1.0908x^4-0.4216x^3 +0.2177x^2+0.1071x$. Here $*$ stands for the global minimizer, and  the Euler's method along the critical point Fingerprint from large $t$ will backward to the true global minimizer.}\label{true-example for six degree}
\end{figure}

\begin{example}
The fingerprint $\mathcal{FP}_1$ of $p(x)=x^6+0.6987x^5-1.0908x^4-0.4216x^3 +0.2177x^2+0.1071x$ illustrated in Fig.\ref{true-example for six degree}
shows that the global minimizer is included in the integral curve to convex $p(x,t)$.
\end{example}

\begin{figure}
  \centering
  \includegraphics[width=9cm]{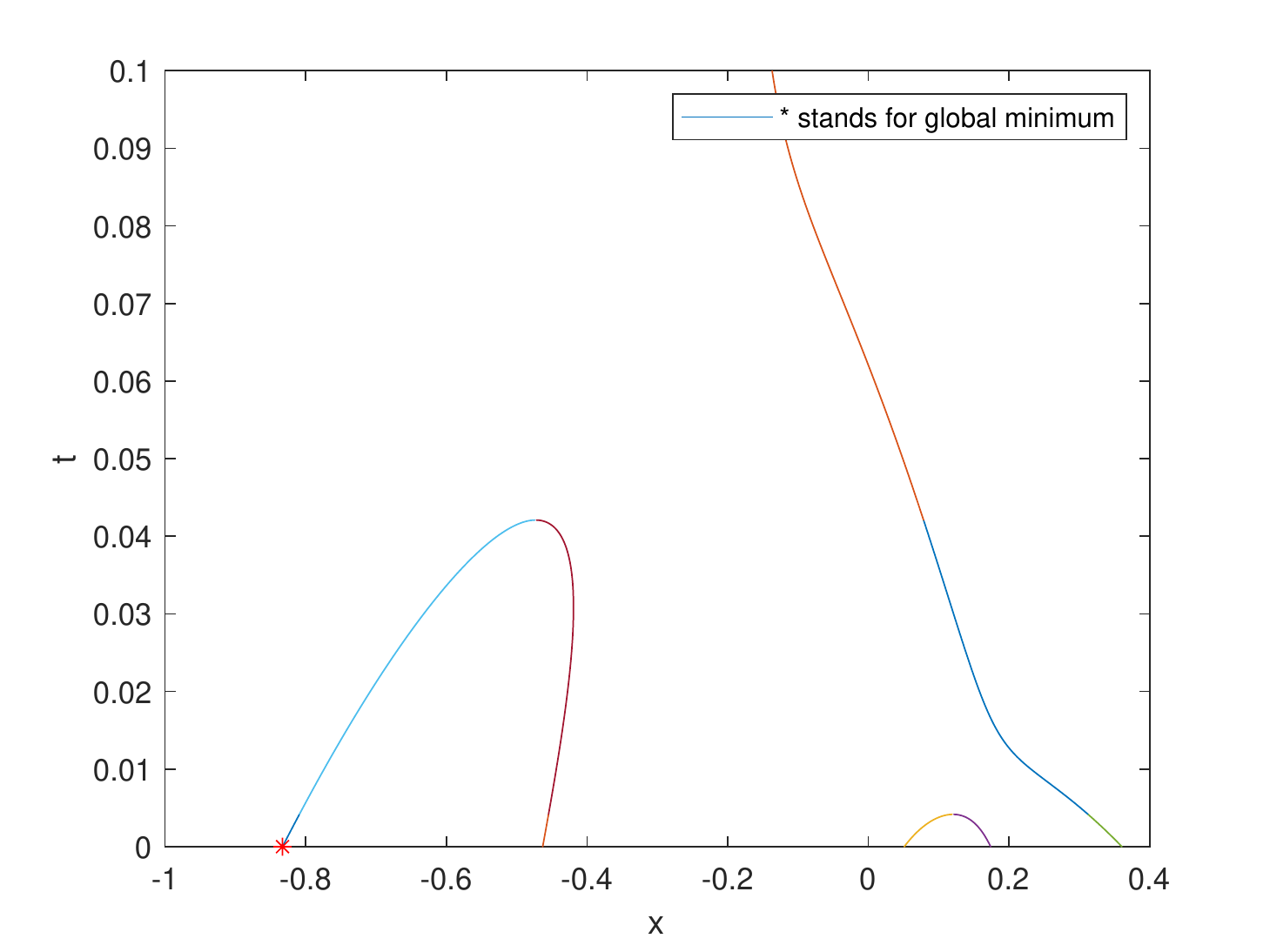}\\
  \caption{A counter-example of fingerprint of $p(x)=x^6-0.8529x^5-0.4243x^4-0.2248x^3 +0.0916x^2-0.0074x$. Here the $*$ stands for the global minimizer, but the most right curve started from large $t>0$, connected only to the local minimizer.}\label{Counter-example for six degree}
\end{figure}

\begin{example}
The fingerprint $\mathcal{FP}_1$ of $p(x)=x^6-0.8529x^5-0.4243x^4-0.2248x^3 +0.0916x^2-0.0074x$ illustrated in Fig.\ref{Counter-example for six degree}
shows that the global minimizer is NOT included in the integral curve to convex $p(x,t)$.
\end{example}

\section{Conclusion}
Motivated by O.Arikan {{\it et al}'s}  \cite{ABK} Stekelov convexification version of backward differential flow algorithm to search the global minimizer of polynomials, we introduce a  natural heat evolution framework, and apply Yuille-Poggio's fingerprint theory as well as their  Yuille-Poggio trajectory equation in computer vision, to characterize the evolution of critical points of $p(x,t)$.  We propose the concepts of {\it escape zone} and {\it confinement zone} of the polynomial $p(x)$, based on which we found  that the global minimizer $x^*$ of a polynomial $p(x)$ can be evolved  inversely  from the global minimizer of its conxification version $p(x,t)=p(x)*g_t(x)$, {\bf{if and only if}} this $x^*$ is in the escape zone of polynomial $p(x)$. 

Specifically, for quartic polynomials, we not only obtained the similar results of O. Arikan et al. \cite{ABK} but also presents a significantly simpler version of Newton's method that can globally minimize quartic polynomials without  convexification.

However, the explicit representation of escape zone and attainable zone of a polynomial $p(x)$ is in general algebraically not tractable, according to the Galois theory.
Thus efficient numerical methods, as well as various criterions of judge the zones are extremely expected.

Moreover, if the global minimizer $x^*$ is outside the escape zone, we have devised a PDE approach to reach this $x^*$, starting from the global minimizer of a para-linear modified version of $p(x)$, which will be disclosed soon. Additionally, findings related to multivariate scenarios will be presented in our upcoming publications.


\begin{landscape}
\begin{figure}
\centering
  \includegraphics[width=22cm]{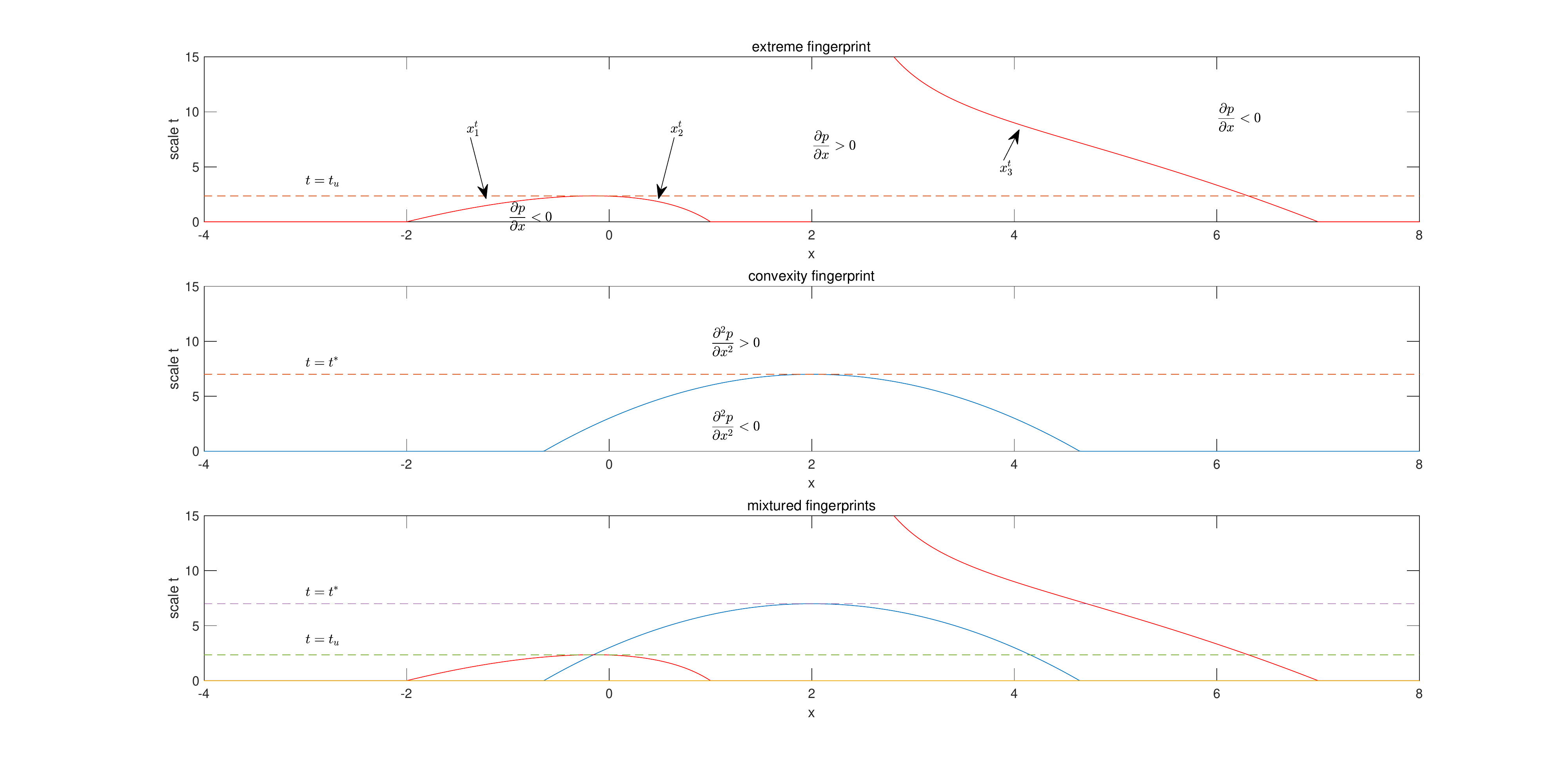}\\
  \caption{Both the fingerprints $\mathcal {FP}_1$ and $\mathcal {FP}_2$ characterize the distribution of critical points and convexity of heat evolved version of a quartic polynomial.}\label{fig:2}
\end{figure}
\end{landscape}

\appendix
\section{Real solutions of cubic equation}\label{appendix-1}

\subsection{Representation by roots}

Recall the classical theory of cubic algebraic equation (cf. \cite{Ze})
\begin{equation}\label{eqn:cubic}
x^3+\alpha x^2+\beta x +\gamma=0,
\end{equation}
According to Newton's method, we have

\begin{lemma}\label{roots-represent-1}
Let $x_i\ (i=1,2,3)$ be the roots (real or complex) of polynomial equation \eqref{eqn:cubic}. We have the following propositions,
\begin{equation}\label{appendix-2}
\begin{split}
  x_1+x_2+x_3&=-\alpha,\\
  x_1x_2+x_2x_3+x_3x_1&=\beta,\\
  x_1x_2x_3&=-\gamma,\\
  x_1^2+x_2^2+x_3^2&=\alpha^2-2\beta,\\
  x_1^3+x_2^3+x_3^3&=-\alpha^3+3\alpha\beta-3\gamma,\\
  x_1^4+x_2^4+x_3^4&=\alpha^4-4\alpha^2\beta+4\alpha^2+2\beta^2.
\end{split}
\end{equation}
\end{lemma}

When the coefficients of the equation are real,  the discriminant of the equation is
\begin{equation}
\Delta=(x_1-x_2)^2(x_2-x_3)^2(x_3-x_1)^2,
\end{equation}
which is equivalent to
\begin{equation}
\Delta=\frac{g^2}{4}+\frac{f^3}{27},
\end{equation}
in which
\begin{equation}\label{eq:d-1}
f=\beta-\frac{\alpha^2}{3}\ \ {\mathrm {and}}\ \
g=\frac{2\alpha^3}{27}-\frac{\alpha\beta}{3}+\gamma.
\end{equation}

\subsection{The real roots described by discriminant}\label{discriminant}
Then the solution of this cubic equation is as follows:

If $\Delta<0$, the equation \eqref{eqn:cubic} contains three distinct real roots.
   \begin{equation}
   \begin{split}
   x_1&=\frac{2}{\sqrt{3}}\sqrt{-f}\sin(\theta)-\frac{\alpha}{3},\\
   x_2&=-\frac{2}{\sqrt{3}}\sqrt{-f}\sin\left (\theta+\frac{\pi}{3} \right )-\frac{\alpha}{3},\\
   x_3&=\frac{2}{\sqrt{3}}\sqrt{-f}\cos\left ( \theta+\frac{\pi}{6} \right )-\frac{\alpha}{3},
   \end{split}
   \end{equation}
   where
   \begin{equation}
   \theta=\frac{1}{3}\arcsin \left ( \frac{3\sqrt{3}g}{2(\sqrt{-f})^3}\right ).
   \end{equation}

  If $\Delta=0$, the solutions contains a single root  and two repeated roots,
  \begin{equation}\label{eqn:cubic-3}
  \begin{split}
  x_1&=-2\left ( \frac{g}{2}\right )^{\frac{1}{3}}-\frac{\alpha}{3},\\
  x_2=x_3&=\left (\frac{g}{2} \right )^{\frac{1}{3}}-\frac{\alpha}{3}.
  \end{split}
  \end{equation}

  Finally, when $\Delta>0$, the equation has only single real root
  \begin{equation}
  x=\left (-\frac{g}{2}+\sqrt{\Delta} \right )^{\frac{1}{3}}+\left (-\frac{g}{2}-\sqrt{\Delta} \right )^{\frac{1}{3}}.
  \end{equation}

\section{Real double roots of quartic polynomials and the structure of $\mathcal{FP}_2\bigcap\mathcal{FP}_3$}\label{solution}
We consider the regularized form of real quartic equation
\begin{equation}\label{App-B-1}
x^4+\beta x^2+\gamma x + \delta =0.
\end{equation}
The structure of its roots is described in \cite{Rees}. It possesses repeated roots if and only if its discriminant $\Delta=0$, where
\begin{equation}
\Delta = 256\gamma^3-128\beta^2\delta^2+144\beta\gamma^2\delta -27\gamma^4+16\beta^4\delta-4\beta^3\gamma^2.
\end{equation}
In addition, we require another four polynomials,
\begin{equation}\label{App-B-2}
\begin{split}
P = & 8\beta,\\
R = & 8\gamma,\\
\Delta_0 = & \beta^2+12\delta,\\
D = & 64\delta-16\beta^2.\\
\end{split}
\end{equation}
The equation \eqref{App-B-1} has one double real roots and two other distinct real roots, if and only if
\begin{equation}\label{App-B-3}
\Delta=0\ \rm{and}\ P<0\ \rm{and}\ D<0\ \rm{and}\ \Delta_0\ne 0.
\end{equation}
The equation \eqref{App-B-1} has one double real roots and a pair of complex roots, if and only if
\begin{equation}\label{App-B-4}
(\Delta = 0\ \rm{and}\ D>0)\ \rm{or}\ (\Delta=0\ \rm{and}\ P>0\ \rm{and}\ (D\ne 0\ \rm{or}\ R\ne 0)).
\end{equation}

To apply the above results to the quartic polynomial evolution,  we may represent the variables according to
\begin{equation}\label{App-B-5}
\begin{split}
\beta = & \frac{2b}{5}+6t,\\
\gamma = & \frac{c}{5},\\
\delta = & \frac{d}{15}+\frac{2b}{5}t+3t^2.\\
\end{split}
\end{equation}

Finally, we see that
\begin{equation}\label{App-B-delta}
\Delta(t)=\sum_{k=0}^6 c_{6-k}t^k,
\end{equation}
in which
\begin{equation}\label{App-B-6-delta}
\begin{split}
c_0=& 27648,\\
 c_1=& \frac{55296b}{5},\\
 c_2 =&\frac{9216b^2}{5},\\
 c_3 =&\frac{4096b^3+1728c^2}{25},\\
 c_4 =& \frac{4864b^4}{625} + \frac{512db^2 + 1728bc^2}{125} - \frac{256d^2}{25},\\
 c_5 = & \frac{32(b^2+5d)(48b^3-80db+135c^2)}{9375},\\
 c_6 = &\frac{256b^4d}{9375} -\frac{32b^3c^2}{3125} - \frac{512b^2d^2}{5625} +\frac{96bc^2d}{625} - \frac{27c^4}{625} +\frac{256d^3}{3375}\\
  \end{split}
  \end{equation}
Here we should explicitly represent  the conditions in \eqref{App-B-3} \eqref{App-B-4}. Actually, based on \eqref{App-B-5}, we have
\begin{equation}
\begin{split}
P< 0&\iff t<-\frac{b}{15},\\
R\ne 0 & \iff c\ne 0,\\
\Delta_0\ne 0&\iff d>\frac{b^2}{5},\\
D > 0\ (=0,\ <0,\ \rm{resp.}) &\iff \left (t+\frac{b}{15}  \right )^2 -\frac{1}{90} \left (d- \frac{b^2}{5} \right )> 0 \ (=0,\  <0,\ \rm{resp.}).
\end{split}
\end{equation}

These implies that
\begin{equation*}
\eqref{App-B-3} \iff -\frac{b}{15}-\frac{1}{\sqrt{90}}\sqrt{d-\frac{b^2}{5}}<t<-\frac{b}{15},\ d-\frac{b^2}{5}>0,\ \rm{and}\ \Delta(t)=0.
\end{equation*}

If we denote
\begin{equation*}
\widehat t=\frac{1}{\sqrt{90}}\sqrt{d-\frac{b^2}{5}}-\frac{b}{15},
\end{equation*}
then we see that
\begin{itemize}
\item If $d-\frac{b^2}{5}<0$, then \eqref{App-B-4} $\iff \Delta(t)=0$.
\item  If $d-\frac{b^2}{5}=0$, then $D>0$ is equivalent to $t\ne -\frac{b}{15}$. Thus

\eqref{App-B-4} $\iff t\ne -\frac{b}{15}$ and $\Delta(t)=0$.

\item If $d-\frac{b^2}{5}>0$, we may classify it into two cases. At first case, if $d\ge \frac{3b^2}{5}$, then $D>0$ is equivalent to
\begin{equation}
0\le t<-\frac{b}{15}- \frac{1}{\sqrt{90}}\sqrt{d-\frac{b^2}{5}},\ \ \rm{or}\ t>-\frac{b}{15}+ \frac{1}{\sqrt{90}}\sqrt{d-\frac{b^2}{5}}.
\end{equation}
Notice that
\begin{equation}
\begin{split}
&D>0\ \lor\ \left ( P>0\ \land\ (D\ne 0\ \lor\ R\ne 0)\right )\\
=&(D>0\ \lor\ P>0)\ \land\ (D>0\ \lor\ D\ne 0\ \lor\ R\ne 0)\\
=&(D>0\ \lor\ P>0)\ \land\ (D\ne 0\ \lor\ R\ne 0)
\end{split}
\end{equation}
which indicates that\newline
\eqref{App-B-4}$\iff \Delta(t)=0, \rm{and}\ t\in \left [0,-\frac{b}{15}-\frac{1}{\sqrt{90}}\sqrt{d-\frac{b^2}{5}}\right )\bigcup \left (-\frac{b}{15},+\infty \right ),$
excluding that both $D=0$ and $c=0$.

But if $\frac{3b^2}{5}> d >\frac{b^2}{5}$, then $D>0$ is equivalent to
\begin{equation}
 t>-\frac{b}{15}+ \frac{1}{\sqrt{90}}\sqrt{d-\frac{b^2}{5}}.
\end{equation}
then
 \eqref{App-B-4} $\iff \Delta(t)=0$  and $t>-\frac{b}{15}$ excluding that both $D=0$ and $c=0$.
\end{itemize}



In a summary, we can obtain at $0\le t_1\le t_2$, respectively corresponds to a dual real roots $x_1$ and $x_2$ of the equation \eqref{App-B-1}, thus
$$\{(x_1,t_1),(x_2,t_2)\}=\mathcal{FP}_2\bigcap \mathcal{FP}_3,$$ in which $t_2\ge t_1$ and $0\le t_1< -\frac{b}{15}$.

%
%
%
%



\begin{thebibliography}{plain}

\bibitem{Las} J. B. Lasserre, Global optimization with polynomials and the problem of moments, SIAM J. Optim. Vol. 11(3), pp:796-817, 2001.




\bibitem{Shor-1} N.Z. Shor, Quadratic optimization problems, Soviet J. Comput. Systems Sci., 25(1987), pp:1-11.

\bibitem{Shor-2} N.Z. Shor, Nondifferentiable optimization and polynomial problems, Kluwer Academic Publishers, 1998

\bibitem{Nef} V.N. Nefedov, Polynomial optimization problem, U.S.S.R. Comput. Maths. Math. Phys., Vol.27, No.3, pp.l3-21, 1987

\bibitem{Zhu-1} J. Zhu and X. Zhang, On global optimizations with polynomials, Optimization Letters, (2008)2: 239-249.

\bibitem{Zhu-2} J. Zhu, S. Zhao and G. Liu, Solution to global minimization of polynomials by backward differential flow,
J. Optim Theory Appl (2014)161: 828-836.

\bibitem{ABK-1} O. Arikan, R.S. Burachik and C.Y. Kaya, "Backward differential flow" may not converge to a global minimizer of polynomials, J. Optim Theory Appl (2015): 167: 401-408.

\bibitem{ABK} O. Arikan, R.S. Burachik and C.Y. Kaya, Steklov regularization and trajectory methods for univariate global optimization, J. Global Optimization, 2019.

\bibitem{BK} Burachik, R.S., Kaya, C.Y. Steklov convexification and a trajectory method for global optimization of multivariate quartic polynomials. Math. Program. (2020). https://doi.org/10.1007/s10107-020-01536-8
\bibitem{Iijima-1} T. Iijima, basic theory of pattern observation, Papers of Tech. Group on Automata and Automatic Control, IEICE, Japan, 1959 (in Japanese).

\bibitem{Iikima-2} T. Iijima, basic theory on normalization of pattern (in case of typical one-dimensional pattern), Bulletin of the Electrotechnical Lab., Vol.26: 368-388, 1962.
\bibitem{Yuille-Poggio-1}     A. L. Yuille and T. Poggio, Fingerprints theorems for zero crossings, J. Opt. Soc. Am. A, 2(5): 683-692, 1985. doi = {10.1364/JOSAA.2.000683}
\bibitem{Yuille-Poggio-2} A. L. Yuille and T. A. Poggio, Scaling theorems for zero crossings, IEEE Trans. on Pattern Analysis and Machine Intelligence, vol.8(1), pp. 15-25, Jan. 1986, doi: 10.1109/TPAMI.1986.4767748.
\bibitem{Ze} E. Zeidler, Oxford users' guide to mathematics, Oxford Univ. Press, 2013.

\bibitem{Kap} Irving Kaplansky, an introduction to differential algebra, Hermann, Paris, 1957.

\bibitem{John} Fritz John, Partial differential equations, Springer-Verlag, 4.ed., 1982.
\bibitem{Walter} W.A. Strauss,  Partial differential equations. John Wiley and Sons Inc. 1992.
\bibitem{Stewart} I.N. Stewart, Galois theory, Chapman and Hall/CRC, 4.ed., 2015.
\bibitem{Rees} E. L. Rees, Graphical discussion of the roots of a quartic equation, The American Mathematical Monthly, 29:2, 51-55, 1922. DOI:~ 10.1080/00029890.1922.11986100
\bibitem{Birk-Mac} Birkhoff, G. and Mac Lane, S. A survey of modern algebra, 3rd ed. New York: Macmillan, 1965.
\bibitem{Nick} Nickalls, R.W.D, The quartic equation: invariants and Euler’s solution revealed, The Mathematical Gazette, vol.93(526), 66-75, 2009.  



\end{thebibliography}
\end{document}